\documentclass[a4paper,12pt,twoside]{article}
\usepackage{amsmath,amssymb,ifthen}
\usepackage[amsmath,thmmarks]{ntheorem}
\usepackage{graphicx}
\usepackage{paper}
\usepackage[all]{xy}
\usepackage{mathrsfs}
\usepackage{url}
\usepackage{mathtools}
\usepackage{stmaryrd}
\usepackage{autobreak}
\usepackage{hhline}
\DeclareMathOperator{\Ind}{Ind}
\DeclareMathOperator{\Res}{Res}

\DeclareMathOperator{\Sp}{Sp}
\DeclareMathOperator{\SO}{SO}
\DeclareMathOperator{\SL}{SL}

\DeclareMathOperator{\lcm}{lcm}
\newcommand{\G}{\mathbb{G}}
\newcommand{\Pcal}{\mathcal{P}}
\newcommand{\Qcal}{\mathcal{Q}}

\SelectTips{cm}{11}
\title{The number of cuspidal representations over a function field and its behavior under base changes}
\author{Takuro Fukayama}

\begin{document}

\maketitle

\begin{abstract}
Let $X$ be a smooth projective curve over a finite field $\F_q$, $k$ be its function field, and $G$ be a simply connected almost simple split group over $\F_q$. We also write $G$ for its structure over $k$. We calculate the sum of multiplicities of all cuspidal representations of $G$ satisfying a given condition assuming the conjectural trace formula. We also observe how the sum changes if we replace $X$ by its base change $X\otimes_{\F_q}\F_{q^m}$.
\end{abstract}

\section{Introduction}
Let $X$ be a smooth projective curve over a finite field $\F_q$, $k$ be its function field, and $G$ be a simply connected almost simple split group over $\F_q$. Given disjoint finite sets $S$, $T$ of places, we are interested in the number of cuspidal representations $\pi$ whose local component $\pi_v$ is the Steinberg representation for $v\in S$, a simple supercuspidal representation for $v\in T$, and unramified elsewhere. More precisely, we want to calculate the sum of multiplicities of such cuspidal representations. 

One of our motivations to calculate the sum of multiplicities is the global Langlands correspondence, which states that there is a correspondence between the set of cuspidal representations of $G$ and the set of $\ell$-adic $\widehat{G}$-local systems on an open subset of $X$. If the considered group is $\GL_n$, then the correspondence is given in \cite{MR1875184}, and roughly speaking, counting cuspidal representations is equivalent to counting $\ell$-adic local systems on an open subset of $\overline{X}=X\otimes_{\F_q}\overline{\F}_q$ which are fixed by the Frobenius $\Fr_q$. Moreover, the base change from $\F_q$ to $\F_{q^m}$ causes the replacement from $\Fr_q$ to $\Fr_{q^m}=\Fr_q^m$, and it is proved in \cite{MR3092473} that the changing law of the number of $\ell$-adic local systems (resp.\ cuspidal representations) can be described as a function of \emph{Lefschetz type}: this is the class of functions $f\colon \Z\to\C$ of the form
\[
f(m)=\sum_{i=1}^r n_i \alpha_i^m
\]
for some $n_1,\ldots, n_r\in\Z$ and $\alpha_1,\ldots,\alpha_r\in\C$. It is also proved in \cite{yu2024elladiclocalsystemshiggs} and \cite{yu2024rank2elladiclocal}  that the number of $\ell$-adic local systems gives a Lefschetz type function under various conditions.  From this perspective, it is expected that the base change from $\F_q$ to $\F_{q^m}$ gives a Lefschetz type function which counts the number of cuspidal representations in general. For example, some results on the number of cuspidal representations of more general algebraic groups are given in \cite{yu2023numbercuspidalautomorphicrepresentations}.

Gross has proposed some conjectures in \cite{MR2843098}:
\begin{conj}[{{\cite[p.~1250]{MR2843098}}}]\label{conj:EP function}
Let $F$ be a non-Archimedean local field, and $f^{EP}$ be the Euler-Poincaré function with respect to a Haar measure $dg$ on $G(F)$. Then 
\[
O_\gamma(f^{EP},dg/dg_\gamma)dg_\gamma=\left\{\begin{array}{lc}
dg_\gamma^{EP} & \text{if}\ \gamma\ \text{is elliptic semisimple}\\
0 & \text{otherwise}
\end{array}\right.
\]
where $dg_\gamma$ is any Haar measure on $G_\gamma(F)$ and $dg_\gamma^{EP}$ is the Euler-Poincaré measure on $G_\gamma(F)$.
\end{conj}

\begin{conj}[{{\cite[Conjecture 5.2]{MR2843098}}}]\label{conj:trace formula}
\sloppy Assume that $\#S\geq 1$ and $\#S+\#T\geq 2$. Then $\varphi=\varphi_{S,T,V}$ has a trace on the discrete spectrum $L_d^2(G)=L_d^2(G(k)\backslash G(\A))$, and this trace is given by a sum of orbital integrals
\[
\Tr(\varphi\mid L_d^2(G))=\sum_{[\gamma]}O_\gamma(\varphi)
\]
where the sum is taken over the $S$-elliptic semisimple torsion conjugacy classes in $G(k)$, and only finitely many terms in the sum are non-zero.
\end{conj}

We will give some explanations. In Conjecture \ref{conj:EP function}, the \emph{Euler-Poincaré function} $f^{EP}$ is defined as follows. Let $\mathcal{B}$ the Bruhat-Tits building of $G(F)$ and $\mathcal{F}$ be the set of facets of $\mathcal{B}$. We write $G(F)_\sigma$ for the stabilizer of each $\sigma\in\mathcal{F}$. Then $G(F)_\sigma$ acts on the set of vertices of $\sigma$ by permutations, and thus one has the sign character $\mathrm{sgn}_\sigma\colon G(F)_\sigma\to \{\pm1\}$. Now we define
\[
f^{EP}\coloneqq\sum_{\sigma\in\mathcal{F}}\frac{(-1)^{\dim\sigma}}{\mathrm{vol}(G(F)_\sigma, dg))}\mathrm{sgn}_\sigma.
\]
Here, $\mathrm{sgn}_\sigma$ is extended to a function on $G(F)$ whose support is $G(F)_\sigma$, which makes $f^{EP}$ into a locally constant function on $G(F)$ with a compact support. The orbital integral of $f^{EP}$ at an element $\gamma\in G(F)$ is defined as
\[
O_\gamma(f^{EP},dg/dg_\gamma)=\int_{G_\gamma(F)\backslash G(F)}f^{EP}(g^{-1}\gamma g)\frac{dg}{dg_\gamma}.
\]
The construction of $f^{EP}$ is due to \cite{MR942522} and Conjecture \ref{conj:EP function} is proved \textit{loc.\ cit.\ }if $F$ is a $p$-adic field. 

In Conjecture \ref{conj:trace formula}, $\varphi=\varphi_{S,T,V}$ is the global test measure on $G(\A)$ constructed in \cite[\S 4]{MR2843098} with respect to finite sets $S$, $T$ of places, and a fixed representation $V=\bigotimes_{v\in S\cup T}V_v$ of the subgroup
\[
G_{S,T}\coloneqq \prod_{v\in S}G(k_v)\times \prod_{v\in T}I_v^+\times \prod_{v\notin S\cup T}G(\mathcal{O}_v)
\]
of $G(\A)$. Here, $\mathcal{O}_v$ is the ring of integers of $k_v$ and $I_v^+$ is a pro-$p$ Iwahori subgroup of $G(\mathcal{O}_v)$. Then $\varphi$ is the product of local measures $\varphi_v$ on $G(k_v)$, and each $\varphi_v$ is constructed so that the traces $\Tr(\varphi_v\mid \pi_v)$ on local representations $\pi_v$ on $G(k_v)$ behave as a characteristic function on the set of such representations which detects a prescribed one. It follows that the trace $\Tr(\phi\mid L_d^2(G))$ counts the number of automorphic representations with given local conditions. (See Section \ref{The sum of multiplicities} for details.) Here Conjecture \ref{conj:trace formula} states that the trace expands as a finite sum of orbital integrals
\[
O_\gamma(\varphi)=\int_{G_\gamma(k)\backslash G(\A)}\varphi(g^{-1}\gamma g)/dz_\gamma
\]
where $dz_\gamma$ is the counting measure on a discrete group $G_\gamma(k)$. We say a semisimple element $\gamma\in G(k)$ is $S$-elliptic if it is elliptic in $G(k_v)$ for all $v\in S$.

If we assume Conjectures \ref{conj:EP function} and \ref{conj:trace formula}, then we can calculate the sum of multiplicities using the Artin $L$-functions determined by $G$, which are defined in Section \ref{section:L-functions}. We consider two cases: 
\begin{enumerate}
\item $S,T\neq\emptyset$,
\item $T=\emptyset$, $\#S\geq 2$.
\end{enumerate}
In the case (1), the sum of multiplicities can be expressed with a single $L$-function $L_{S,T}(M_G)$, and we will see that it gives a Lefschetz type function by considering base changes. For the precise statement and proof, see Proposition \ref{prop:the sum of multiplicities gives a Lefschetz type function}. On the other hand, it is expressed as a sum of finitely many $L$-functions $L_{S,T}(M_{G_\gamma})$ in the case (2), and the situation becomes more complicated. Our main work in this paper is to analyze that sum of $L$-functions for $G=\SL_n$ and $\Sp_{2n}$. Moreover, we proved that the sum gives Lefschetz type functions under base changes, in some special cases. Here are our main results:

\begin{thm}[Theorem \ref{thm:the sum of L-functions for SL_l}]
Let $\ell$ be a prime and $G=\SL_\ell$. Then the sum $\sum_{[\gamma]}L_S(M_{G_{\gamma}})$ gives a Lefschetz type function in $m$ under the base change from $\F_q$ to $\F_{q^m}$.
\end{thm}

\begin{thm}[Theorems \ref{thm:the sum of L-functions for Sp_4} and \ref{thm:the sum of L-functions for Sp_6}]
Let $G=\Sp_4$ or $\Sp_6$. Then the sum $\sum_{[\gamma]}L_S(M_{G_{\gamma}})$ gives a Lefschetz type function in $m$ under the base change from $\F_q$ to $\F_{q^m}$.
\end{thm}

The outline of this paper is as follows. In Section \ref{Lefschetz type functions}, we will see some examples of Lefschetz type functions. In Section \ref{section:L-functions}, we review the notion of $L$-functions of motives for algebraic groups and observe their behaviors under base changes. In Section \ref{The sum of multiplicities}, we explain how to calculate the sum of multiplicities using the trace formula, following the argument in \cite{MR2843098}. We also see that the behavior of the sum under base changes gives a Lefschetz type function in $S,T\neq\emptyset$ case. In Section \ref{The sum of L-functions: SL_n case}, we calculate the $L$-functions for semisimple centralizers for $G=\SL_n$ and see that the sum of them gives a Lefschetz type function under base changes if $n=\ell$ is prime. We also give another result related to the sum of $L$-functions. In Section \ref{The sum of L-functions: Sp_2n case}, we calculate the $L$-functions for semisimple centralizers for $G=\Sp_{2n}$ and see that the sum of them gives a Lefschetz type function under base changes if $n=2,3$.

\subsection*{Notation}
\begin{enumerate}
\item Let $\F_q$ be a finite field of characteristic $p$ with $q$ elements, $X$ be a smooth projective curve over $\F_q$, $k$ be its function field, and $\A=\A_k$ be an adele ring of $k$. For each place $v$ of $k$, let $k_v$ be the corresponding complete discrete valuation field with the ring of integers $\mathcal{O}_v$ and the residue field $\kappa_v$.  We fix disjoint finite sets $S$, $T$ of places of $k$ and assume that $S$ is always non-empty while $T$ can be empty.  

\item Considering the base change from $\F_q$ to $\F_{q^m}$ for $m\in\Z_{>0}$, we put $X_m\coloneqq X\otimes_{\F_q}\F_{q^m}$, $k_m\coloneqq k\otimes_{\F_q}\F_{q^m}$, $S_m\coloneqq S\otimes_{\F_q}\F_{q^m}$, $T_m\coloneqq T\otimes_{\F_q}\F_{q^m}$ so that $X=X_1$, $k=k_1$, $S=S_1$, $T=T_1$. We also put $\overline{X}\coloneqq X\otimes_{\F_q}\overline{\F}_q$, $\overline{S}\coloneqq S\otimes_{\F_q}\overline{\F}_q$, $\overline{T}\coloneqq T\otimes_{\F_q}\overline{\F}_q$.

\item Let $G$ be a connected reductive group over $\F_q$, and assume that $G$ is simply connected almost simple split except for Section \ref{section:L-functions}. We also write $G$ for its structure over $k$. In the second half of this paper, we mainly consider the cases $G=\SL_n$, $\Sp_{2n}$. We write $G_\gamma$ for the centralizer of an element $\gamma\in G$. 
\end{enumerate}

\subsection*{Acknowledgement}
I would like to thank my supervisor, Professor Yoichi Mieda, for providing useful advice to support my study.

\section{Lefschetz type functions}\label{Lefschetz type functions}
Our aim in this paper is to clarify the way the number of cuspidal representations with prescribed local conditions changes under base changes, which is expected to follow the law called Lefschetz type. This notion arises from the Grothendieck-Lefschetz trace formula, which gives a formula to count the number of $\F_{q^m}$-rational points of a given variety over $\F_q$. In this section, we introduce the definition and some examples of Lefschetz type functions.
\begin{defn}
A function $f\colon \Z_{>0}\to \C$ is called \emph{of Lefschetz type} if it is of the form
\[
f(m)=\sum_{i=1}^r n_i\alpha_i^m
\]
for some $n_1,\ldots, n_r\in\Z$ and $\alpha_1,\ldots,\alpha_r\in\C$. 
\end{defn}

\begin{exa}\label{exa:Lefschetz type function}
	\begin{enumerate}
	\item If $f_1$ and $f_2$ are Lefschetz type functions, then so are $f_1\pm f_2$, $f_1\cdot f_2$.
	\item Let $n$ be a positive integer and define the function $\chi_n$ by $\chi_n(m)=n$ for $n\mid m$ and $\chi_n(m)=0$ for $n\nmid m$. Then $\chi_n$ is of Lefschetz type. Indeed, it has the expression $\chi_n(m)=\sum_{i=0}^{n-1}\zeta_n^{i m}$ for $\zeta_n=e^{2\pi\sqrt{-1}/n}$.
	\item Let $Y$ be an algebraic variety over $\F_q$. Then the function $m\mapsto \#Y(\F_{q^m})$ is of Lefschetz type by the Grothendieck-Lefschetz trace formula.
	\end{enumerate}
\end{exa}

\begin{prop}\label{prop:Lefschetz type function}
For a positive integer $N$ and a Lefschetz type function $f$, define a function $f_N$ by $f_N(m)\coloneqq f(\lcm(N, m))^{\gcd(N, m)}$. Then $f_N$ is of Lefschetz type.
\end{prop}

\begin{prf}
If $N=1$, the assertion is trivial. First we consider the case where $N=\ell^e$ for some prime $\ell$ and $e\in\Z_{>0}$. We have 
\[
f_{\ell^e}(m)=\left\{\begin{array}{ll}
f(\ell^{e-v_\ell(m)}m)^{\ell^{v_\ell(m)}} & \text{if}\ v_\ell(m)<e\\
f(m)^{\ell^e} & \text{if}\ v_\ell(m)\geq e
\end{array}\right.
\]
where $v_\ell(m)$ is the valuation of $m$ at $\ell$. Define another Lefschetz type function by $g(m)\coloneqq f(\ell m)$. Then we have
\[
g_{\ell^{e-1}}(m)=\left\{\begin{array}{ll}
g(\ell^{e-1-v_\ell(m)} m)^{\ell^{v_\ell(m)}}=f(\ell^{e-v_\ell(m)}m)^{\ell^{v_\ell(m)}} & \text{if}\ v_\ell(m)<e-1\\
g(m)^{\ell^{e-1}}=f(\ell m)^{\ell^{e-1}} & \text{if}\ v_\ell(m)\geq e-1.
\end{array}\right.
\]
It follows that
\[
f_{\ell^e}(m)-g_{\ell^{e-1}}(m)=\left\{\begin{array}{ll}
0 & \text{if}\ \ell^e\nmid m\\
f(m)^{\ell^e}-f(\ell m)^{\ell^{e-1}} & \text{if}\ \ell^e\mid m.
\end{array}\right.
\]
Since $g_{\ell^{e-1}}$ is of Lefschetz type by induction hypothesis, it suffices to show that so is $f_{\ell^e}-g_{\ell^{e-1}}$. Write as $f(m)=\sum_i n_i\alpha_i^m$ and regard each $\alpha_i$ as an indeterminate. Then one has $f(m)^\ell\equiv\sum_i n_i^\ell\alpha_i^{\ell m}\equiv\sum_i n_i\alpha_i^{\ell m}=f(\ell m)\bmod\ell$. We claim that $f(m)^{\ell^e}\equiv f(\ell m)^{\ell^{e-1}}\bmod\ell^e$. Indeed, if we put $s\coloneqq f(m)^{\ell^{e-1}}$, $t\coloneqq f(\ell m)^{\ell^{e-2}}$, then $s\equiv t\bmod\ell^{e-1}$ by induction hypothesis and we have
\[
s^{\ell-1}+s^{\ell-2}t+\cdots+t^{\ell-1}\equiv  s^{\ell-1}+s^{\ell-1}+\cdots+s^{\ell-1}=\ell s^{\ell-1}\equiv 0\bmod\ell
\]
which implies
\[
s^\ell-t^\ell=(s-t)(s^{\ell-1}+s^{\ell-2}t+\cdots+t^{\ell-1})\equiv 0\bmod \ell^e.
\]
\sloppy This means that there exists a Lefschetz type function $h$ such that $f(m)^{\ell^e}-f(\ell m)^{\ell^{e-1}}=\ell^e\cdot h(m)$ and so that we have
\[
f_{\ell^e}-g_{\ell^{e-1}}=\chi_{\ell^e}\cdot h
\]
for $\chi_{\ell^e}$ defined in Example \ref{exa:Lefschetz type function}. So we are done.

Now we consider the general case. Suppose $f_N$ is of Lefschetz type. It suffices to show that $f_{\ell^eN}$ is also of Lefschetz type for any prime $\ell$ not dividing $N$, and for any $e\in\Z_{>0}$. One can see that
\begin{align*}
\gcd(\ell^eN,m)=\gcd(N,m)\gcd(\ell^e,m),\ \ \ 
\lcm(\ell^eN,m)=\lcm(N,\lcm(\ell^e,m)).
\end{align*}
It follows that 
\begin{align*}
f_{\ell^eN}(m)&=f(\lcm(\ell^eN,m))^{\gcd(\ell^eN,m)}\\
&=f(\lcm(N,\lcm(\ell^e,m)))^{\gcd(N,m)\gcd(\ell^e,m)}\\
&=f_N(\lcm(\ell^e,m))^{\gcd(\ell^e,m)}\\
&=(f_N)_{\ell^e}(m).
\end{align*}
We have already shown that $(f_N)_{\ell^e}$ is of Lefschetz type if so is $f_N$.
\end{prf}

\begin{cor}\label{cor:Lefschetz type function}
Let $k$ be a global function field over $\F_q$, $T$ be a finite set of places, and $f$ be a Lefschetz type function. Then the function $m\mapsto \prod_{w\in T_m}f(m\deg w)$ is also of Lefschetz type.
\end{cor}

\begin{prf}
For each $v\in T$, we have the following Cartesian diagram:
\[
\xymatrix{
\{w_1,\ldots,w_r\}\ar[r]\ar[d]& T_m\ar[d]\\
v\ar[r] & T\lefteqn{.}
}
\]
Here, $w_1,\ldots, w_r$ are distinct places of $k_m$ lying over $v$ with $r=\gcd(\deg v, m)$, $\deg w_i=\deg v/r$. It follows that
\begin{align*}
\prod_{w\in T_m}f(m\deg w)&=\prod_{v\in T}\prod_{\substack{w\in T_m\\w\mid v}}f(m\deg w)\\
&=\prod_{v\in T}f(m\deg v/r)^r\\
&=\prod_{v\in T}f(\lcm(\deg v,m))^{\gcd(\deg v,m)}.
\end{align*}
The assertion follows from Proposition \ref{prop:Lefschetz type function}.
\end{prf}

\section{$L$-functions of motives}\label{section:L-functions}
In this section, we assume that $G$ is a connected reductive group over $\F_q$. We are going to explain the motive for $G$ and associated $L$-functions. See \cite{MR1474159} for general definitions.

Let $S_0$ be a maximal $\F_q$-split torus of $G$, and $T_0$ be its centralizer in $G$. Then $\Gamma_q\coloneqq\Gal(\overline{\F}_q/\F_q)$ acts on the character group $X^*(T_0)\coloneqq\Hom_{\overline{\F}_q}(T_0,\mathbb{G}_m)$. Moreover, $V\coloneqq X^*(T_0)\otimes\Q$ has the structure $V=\bigoplus_{d\geq 1}V_d$ of a graded $\Gamma_q$-representation over $\Q$ in a canonical way \cite[p.~289]{MR1474159}. Now we define the Artin-Tate motive $M_G$ for $G$ to be
\[
M_G\coloneqq\bigoplus_{d\geq1}V_d(1-d).
\]
Here, each $V_d(1-d)=V_d\otimes\Q(-1)^{\otimes d-1}$ is an Artin-Tate motive of weight $2(d-1)$. 

\begin{lem}[{{\cite[Lemma 2.1]{MR1474159}}}]\label{lem:some properties of motives}
	\begin{enumerate}
	\item If $G$ is isogenous to $G'$ over $\F_q$, then $M_G=M_{G'}$.
	\item $M_{G_1\times G_2}=M_{G_1}\oplus M_{G_2}$.
	\item Let $K/\F_q$ be a finite extension and $G_K$ a connected reductive group over $K$. Then $M_{\Res_{K/\F_q}G_K}= \Ind_{K/\F_q}M_{G_K} $.
 	\end{enumerate}
\end{lem}

\begin{exa}\label{exa:motives}
	\begin{enumerate}
	\item $M_{\GL_n}=\Q\oplus\Q(-1)\oplus\Q(-2)\oplus\cdots\oplus\Q(1-n)$.
	\item $M_{\SL_n}=\Q(-1)\oplus\Q(-2)\oplus\Q(-3)\oplus\cdots\oplus\Q(1-n)$.
	\item $M_{\U_n}=\Q[\sigma]\oplus\Q(-1)\oplus\Q[\sigma](-2)\oplus\cdots\oplus\Q[\sigma^n](1-n)$. Here, $\sigma$ is the non-trivial quadratic character of $\Gamma_q$, namely, $\Q[\sigma]$ is the sign representation.
	\item $M_{\Sp_{2n}}=\Q(-1)\oplus\Q(-3)\oplus\Q(-5)\oplus\cdots\oplus\Q(1-2n)$.
	\item $M_{\SO_{2n+1}}=\Q(-1)\oplus\Q(-3)\oplus\Q(-5)\oplus\cdots\oplus\Q(1-2n)$.
	\end{enumerate}
	
(1), (4), (5) are given in \cite[\S 2]{MR1474159}. For (2), the canonical isogeny $\SL_n\times\G_m\to\GL_n$ gives the identity $M_{\SL_n}\oplus \Q=M_{\GL_n}$ by Lemma \ref{lem:some properties of motives}.
	
Now we explain (3). Let $V=\bigoplus_{d\geq 1}V_d$ be the graded $\Gamma_q$-representation attached to $\U_n$. Then it corresponds to the one attached to $\GL_n$ as graded vector spaces because we have $\U_{n,\overline{\F}_q}\cong \GL_{n,\overline{\F}_q}$. Thus $V$ can be regarded as a quotient of the space of symmetric polynomials over $\Q$ in $n$ variables of positive degree, according to the construction in \cite[p.~289]{MR1474159}. Since the Frobenius $\Fr_q$ acts on it by multiplying each variable by $-1$, it follows that $V_d=\Q[\sigma]$ if $d$ is odd and $V_d=\Q$ if $d$ is even. 
\end{exa}

Now we define the $L$-functions of motives. Let $X$ be a smooth projective curve over $\F_q$, $k$ be its function field, and $S$, $T$ be disjoint finite sets of places. For any place $v$, each $V_d$ can be regarded as a representation of $\Gamma_v\coloneqq\Gal(\overline{\kappa_v}/\kappa_v)$ via the inclusion $\Gamma_v=\Gal(\overline{\F}_q/\F_{q^{\deg v}})\subset\Gamma_q$. Let $\Fr_v=\Fr_{q^{\deg v}}\in\Gamma_v$ be the geometric Frobenius element.

\begin{defn}
	\begin{enumerate}
	\item For each place $v$, define the local $L$-function of $M_G$ as
	\begin{align*}
	L_v(M_{G},s)&\coloneqq\det(1-q^{-s\deg v}\Fr_v\mid M_G)^{-1}\\
	&\coloneqq\prod_{d\geq 1}\det(1-q^{-(s+1-d)\deg v}\Fr_v\mid V_d)^{-1}
	\end{align*}
	for $s\in\C$, and the global $L$-function of $M_G$ by
	\[
	L(M_{G},s)\coloneqq\prod_vL_v(M_{G},s)
	\]
	for $\mathrm{Re}\,s\gg 0$. It can be extended to the meromorphic function on $\C$.
	\item Define the $(S,T)$-modified $L$-function of $M_G$ as
	\[
	L_{S,T}(M_{G},s)\coloneqq\frac{L(M_G,s)}{\prod_{v\in S}L_v(M_G,s)\prod_{v\in T}L_v(M_G,s-1)}.
	\]
	We denote its value at $s=0$ by $L_{S,T}(M_G)$.
	\end{enumerate}
\end{defn}

\begin{rem}\label{rem:L-function}
Let $H^i(X)$, $H^i(S)$, $H^i(T)$ be the $\ell$-adic cohomologies of $\overline{X}$, $\overline{S}$, $\overline{T}$ for $\ell$ not dividing $q$. Then $H^i(X)\otimes M_G$, $H^i(S)\otimes M_G$, $H^i(T)\otimes M_G$ become $\ell$-adic motives and we have the explicit formula:
\begin{align*}
L_{S,T}(M_{G},s)&=\det(1-q^{-s}\Fr_q\mid H^1(X)\otimes M_G)\\
&\quad\quad\cdot\frac{\det(1-q^{-s}\Fr_q\mid H^0(S)\otimes M_G)}{\det(1-q^{-s}\Fr_q\mid M_G)}\\
&\quad\quad\cdot\frac{\det(1-q^{1-s}\Fr_q\mid H^0(T)\otimes M_G)}{\det(1-q^{1-s}\Fr_q\mid M_G)}.
\end{align*}
This follows from the expression
\[
L(M_G,s)=\frac{\det(1-q^{-s}\Fr_q\mid H^1(X)\otimes M_G)}{\det(1-q^{-s}\Fr_q\mid M_G)\det(1-q^{1-s}\Fr_q\mid M_G)}
\]
in \cite[p.~1250]{MR2843098} and the identities
\begin{align*}
\det(1-q^{-s}\Fr_q\mid H^0(S)\otimes M_G)&=\prod_{v\in S} \det(1-q^{-s\deg v}\Fr_v\mid M_G)\\
\det(1-q^{1-s}\Fr_q\mid H^0(T)\otimes M_G)&=\prod_{v\in T} \det(1-q^{(1-s)\deg v}\Fr_v\mid M_G).
\end{align*}
If $S$ (resp.\ $T$) is non-empty, then $\Fr_q$ on $H^0(S)$ (resp.\ $H^0(T)$) has at least one $1$ for its eigenvalue. Indeed, $\Fr_q$ acts on it as a permutation matrix. It follows that the second (resp.\ third) factor in the above formula is holomorphic on $\C$.
\end{rem}

\begin{exa}
Consider the case where $X=\P^1_{\F_q}$ and $S=\{0\}$, $T=\{\infty\}$ with $\deg 0=\deg \infty=1$. Then $H^1(X)=0$ and $H^0(S)\simeq H^0(T)\simeq\Q_\ell$ are the trivial representations, and thus $L_{S,T}(M_G)=1$ by the formula in Remark \ref{rem:L-function}
\end{exa}

Define the rational function $Z_q(t)=\prod_{d\geq1}Z_{q,d}(t)$ by 
\begin{align*}
Z_{q,d}(t)&=\det(1-t^{d-1}\Fr_q\mid H^1(X)\otimes V_d)\\
&\quad\quad\cdot\frac{\det(1-t^{d-1}\Fr_q\mid H^0(S)\otimes V_d)}{\det(1-t^{d-1}\Fr_q\mid V_d)}\\
&\quad\quad\cdot\frac{\det(1-t^d\Fr_q\mid H^0(T)\otimes V_d)}{\det(1-t^d\Fr_q\mid V_d)}
\end{align*}
so that $Z_q(q)=L_{S,T}(M_G)$. 

\begin{prop}\label{prop:L-functions of motives}
Assume that $S$ and $T$ are non-empty. Then $Z_q(t)$ lies in $\Z[t]$. Moreover, if we consider the base change from $\F_q$ to $\F_{q^m}$ and denote the motive of $G$ over $\F_{q^m}$ by $M_{G, m}$, then the function $m\mapsto L_{S_m,T_m}(M_{G, m})$ is of Lefschetz type.
\end{prop}

\begin{prf}
Let $J_X$ (resp.\ $J_S$, $J_T$) be the multiset of eigenvalues of $\Fr_q$ on $H^1(X)$ (resp.\ $H^0(S)$, $H^0(T)$). Note that $J_X,J_S,J_T\subset\overline{\Z}$. We have
\begin{align*}
Z_{q,d}(t)&=\prod_{\alpha\in J_X}\det(1-\alpha t^{d-1}\Fr_q\mid V_d)\\
&\quad\quad \cdot\frac{\prod_{\beta\in J_S}\det(1-\beta t^{d-1}\Fr_q\mid V_d)}{\det(1-t^{d-1}\Fr_q\mid V_d)}\cdot\frac{\prod_{\gamma\in J_T}\det(1-\gamma t^d\Fr_q\mid V_d)}{\det(1-t^d\Fr_q\mid V_d)}\\
&=\prod_{\alpha\in J_X}\det(1-\alpha t^{d-1}\Fr_q\mid V_d)\\
&\quad\quad \cdot\prod_{\beta\in J_S-\{1\}}\det(1-\beta t^{d-1}\Fr_q\mid V_d)\cdot\prod_{\gamma\in J_T-\{1\}}\det(1-\gamma t^d\Fr_q\mid V_d).
\end{align*}
Let $J_d$ be the multiset of eigenvalues of $\Fr_q$ on each $V_d$. Then each $J_d$ consists of roots of unity, and one has
\[
\det(1-u \Fr_q\mid V_d)=\prod_{\zeta\in J_d}(1-\zeta u)
\]
for an indeterminate $u$. Thus putting $J\coloneqq (J_X\cup J_S\cup J_T-\{1,1\})\cup\bigcup_{d\geq1}J_d$, the above formula shows that $Z_q(t)=\prod_{d\geq1} Z_{q,d}(t)$ lies in $\Z[t,J]$. Since $Z_q(t)$ is fixed by $\Gal(\overline{\Q}/\Q)$, it is a polynomial in $t$ over $\Z$. 

Consider the base change from $\F_q$ to $\F_{q^m}$. The above argument shows that $Z_q(t)$ is a polynomial in $\{t\}\cup J$ over $\Z$, so we write $Z_q(t,J)$ for it. One can see that $M_{G, m}$ (resp.\ $H^1(X_m)$, $H^0(S_m)$, $H^0(T_m)$) has the same underlying vector space as $M_{G}$ (resp.\ $H^1(X)$, $H^0(S)$, $H^0(T)$), and the action of $\Fr_{q^m}$ on the former corresponds to the action of $\Fr_q^m$ on the latter. It follows that $J_{X_m}=J_X^m$, $J_{S_m}=J_S^m$, $J_{T_m}=J_T^m$, and that $L_{S_m,T_m}(M_{G,m})=Z_q(q^m,J^m)$. Thus our assertion follows.
\end{prf}

\section{The sum of multiplicities}\label{The sum of multiplicities}
Let $G$ be a simply connected almost simple split group over a function field $k$ of characteristic $p$. We will consider the sum of multiplicities of cuspidal representations satisfying a prescribed local behavior. The argument follows mainly from \cite{MR2843098}. 

Fix a family $\chi=\{\chi_v\}_{v\in T}$ of affine generic characters (see \cite[9.2]{MR2730575} for the definition) $\chi_v\colon I_v^+\to\C^\times$, where $I_v^+$ is a pro-$p$ Iwahori subgroup of $G(k_v)$. Let $\mathrm{Irr}_{0}^{S,T}(G,\chi)$ denote the set of isomorphic classes of cuspidal representations $\pi$ of $G$ satisfying the following conditions:
\begin{itemize}
\item[\textbullet] For $v\in S$, $\pi_v$ is isomorphic to the Steinberg representation.
\item[\textbullet] For $v\in T$, $\pi_v$ is isomorphic to an irreducible component of $\Ind_{I_v^+}^{G(k_v)}\chi_v$. (In particular, $\pi_v$ is simple supercuspidal.)
\item[\textbullet] For $v\notin S\cup T$, $\pi_v$ is unramified.
\end{itemize}
We also consider the set $\mathrm{Irr}_0^{S,T}(G)\coloneqq\bigcup_\chi\mathrm{Irr}_{0}^{S,T}(G,\chi)$ of isomorphic classes of cuspidal representations $\pi$ of $G$ satisfying the conditions:
\begin{itemize}
\item[\textbullet] For $v\in S$, $\pi_v$ is isomorphic to the Steinberg representation.
\item[\textbullet] For $v\in T$, $\pi_v$ is a simple supercuspidal representation.
\item[\textbullet] For $v\notin S\cup T$, $\pi_v$ is unramified.
\end{itemize}
If $T=\emptyset$, then $\mathrm{Irr}_{0}^{S,\emptyset}(G,\chi)=\mathrm{Irr}_0^{S,\emptyset}(G)$ by definition and we just write $\mathrm{Irr}_0^{S}(G)$ for it.

Now we are going to calculate $\sum_{\pi\in\mathrm{Irr}_{0}^{S,T}(G,\chi)}m(\pi)$ and $\sum_{\pi\in\mathrm{Irr}_0^{S,T}(G)}m(\pi)$. We will use the global test measure $\varphi=\varphi_{S,T,V}$ on $G(\A)$ defined in \cite[p.~1241]{MR2843098}. In our case, $V$ corresponds to the family $\chi=\{\chi_v\}_{v\in T}$. 

\begin{prop}[{{\cite[pp.~1246--1249]{MR2843098}}}]\label{prop:the sum of multiplicities}
	\begin{enumerate}
	\item Assume that $\#S\geq 1$ and\\ $\#T\geq 1$. Then we have
	\[
	\Tr(\varphi\mid L_d^2(G))=\Tr(\varphi\mid L_0^2(G))=\sum_{\pi\in\mathrm{Irr}_{0}^{S,T}(G,\chi)}m(\pi).
	\]
	\item Assume that $\#S\geq 2$ and $T=\emptyset$. Then we have
	\[
	\Tr(\varphi\mid L_d^2(G))=1+(-1)^{\#S\cdot r(G)}\Tr(\varphi\mid L_0^2(G))=1+(-1)^{\#S\cdot r(G)}\sum_{\pi\in\mathrm{Irr}_{0}^{S}(G)}m(\pi).
	\]
	Here, $r(G)$ stands for the rank of $G$.
	\end{enumerate}
\end{prop}

By Proposition \ref{prop:the sum of multiplicities}, one can calculate the sum of multiplicities using $L$-functions defined in Section \ref{section:L-functions} if we assume Conjectures \ref{conj:EP function} and \ref{conj:trace formula}. When $T$ is empty, we will write simply $L_S(-)$ for $L_{S,\emptyset}(-)$.

\begin{thm}[{{\cite[Theorem 7.1, p.~1251]{MR2843098}}}]\label{thm:the sum of multiplicities} \sloppy Assume that Conjectures \ref{conj:EP function} and \ref{conj:trace formula} are true. Then we have the following:
	\begin{enumerate}
	\item Assume that $\#S\geq 1$ and $\#T\geq 1$. Then we have
	\[
	\sum_{\pi\in\mathrm{Irr}_{0}^{S,T}(G,\chi)}m(\pi)=(-1)^{(\#S+\#T)\cdot r(G)}L_{S,T}(M_G).
	\]
	\item Assume that $\#S\geq 2$ and $T=\emptyset$. Then we have
	\[
	1+(-1)^{\#S\cdot r(G)}\sum_{\pi\in\mathrm{Irr}_{0}^{S}(G)}m(\pi)=\sum_{[\gamma]}L_S(M_{G_\gamma})
	\]
	where the sum on the right-hand side is taken over the semisimple conjugacy classes in $G(\F_q)$. 
	\end{enumerate}
\end{thm}

\begin{rem}\label{rem:the number of simple supercuspidal representations}
\sloppy If Conjectures \ref{conj:EP function} and \ref{conj:trace formula} are true, then the sum of multiplicities $\displaystyle\sum_{\pi\in\mathrm{Irr}_{0}^{S,T}(G,\chi)}m(\pi)$ does not depend on the choice of $\chi=\{\chi_v\}_{v\in T}$ by Theorem \ref{thm:the sum of multiplicities}. In our case, $\Ind_{I_v^+}^{G(k_v)}\chi_v$ decomposes into the direct sum of $\#Z(\F_{q^{\deg v}})$ distinct simple supercuspidal representations. Here, $Z$ is the center of $G$, which is finite by our assumption. On the other hand, there are $(\#Z(\F_{q^{\deg v}}))^2(q^{\deg v}-1)$ isomorphic classes of simple supercuspidal representations of $G(k_v)$ \cite[p.~63]{MR2730575}. So we have 
\[
\sum_{\pi\in\mathrm{Irr}_0^{S,T}(G)}m(\pi)=\left(\prod_{v\in T}\#Z(\F_{q^{\deg v}})\cdot(q^{\deg v}-1)\right)\cdot \sum_{\pi\in\mathrm{Irr}_{0}^{S,T}(G,\chi)}m(\pi).
\]
\end{rem}

Now we observe how the sum of multiplicities changes under the base change from $\F_q$ to $\F_{q^m}$. First, we consider the case $S,T\neq\emptyset$. Then we have the following result:

\begin{prop}\label{prop:the sum of multiplicities gives a Lefschetz type function}
Assume that Conjectures \ref{conj:EP function} and \ref{conj:trace formula} are true, and that $S$ and $T$ are non-empty. Then the function
\[
m\mapsto\sum_{\pi\in\mathrm{Irr}_{0}^{S_m,T_m}(G_{k_m})}m(\pi)
\]
is of Lefschetz type.
\end{prop}

\begin{prf}
Due to Theorem \ref{thm:the sum of multiplicities} (1) and Remark \ref{rem:the number of simple supercuspidal representations}, we have
\begin{align*}
&\sum_{\pi\in\mathrm{Irr}_{0}^{S_m,T_m}(G_{k_m})}m(\pi)\\
&\quad=\left(\prod_{w\in T_m}\#Z(\F_{(q^m)^{\deg w}})((q^m)^{\deg w}-1)\right)\cdot\left(\prod_{v\in S_m\cup T_m}(-1)^{r(G)}\right)\cdot L_{S_m,T_m}(M_{G,m}).
\end{align*}
The functions
\begin{align*}
&m\mapsto\prod_{w\in T_m}\#Z(\F_{q^{m\deg w}})(q^{m\deg w}-1)\\
&m\mapsto\prod_{v\in S_m\cup T_m}(-1)^{r(G)}
\end{align*}
are of Lefschetz type by Corollary \ref{cor:Lefschetz type function}. Since
\[
m\mapsto L_{S_m,T_m}(M_{G,m})
\]
is also a Lefschetz type function by Proposition \ref{prop:L-functions of motives}, the assertion follows.
\end{prf}

In the rest of the paper, we consider the case $\#S\geq2$ and $T=\emptyset$, where it is expected that the sum of multiplicities can be calculated by the sum
\[
\sum_{[\gamma]}L_S(M_{G_\gamma})
\]
 of $L$-functions in Theorem \ref{thm:the sum of multiplicities}. We continue to assume that $G$ is simply connected and almost simple. For each $\gamma$, define the rational function $Z_{q,\gamma}(t)=\prod_{d\geq 1}Z_{q,\gamma, d}(t)$ by
\begin{align*}
Z_{q,\gamma, d}(t)=\det(1-t^{d-1}\Fr_q\mid H^1(X)\otimes V_{\gamma, d})\cdot\frac{\det(1-t^{d-1}\Fr_q\mid H^0(S)\otimes V_{\gamma, d})}{\det(1-t^{d-1}\Fr_q\mid V_{\gamma, d})}
\end{align*}
where $M_{G_\gamma}=\bigoplus_{d\geq 1}V_{\gamma, d}(1-d)$. Then we have $Z_{q,\gamma, d}(q)=L_S(M_{G_\gamma})$. Note that  each $Z_{q,\gamma}(t)$ is no longer a polynomial but a rational function, unlike the case where $S$ and $T$ are non-empty. More concretely, we can write as
 \begin{align*}
 L_S(M_{G_\gamma})&=\frac{\det(1-\Fr_q\mid M_{G_\gamma})}{\det(1-q\Fr_q\mid M_{G_\gamma})}\cdot\det(1-\Fr_q\mid H^1(X)\otimes M_{G_{\gamma}})\\
 &\quad\quad\cdot\frac{\det(1-\Fr_q\mid H^0(S)\otimes M_{G_{\gamma}})}{\det(1-\Fr_q\mid M_{G_\gamma})^2}\\
 &=\frac{\det(1-\Fr_q\mid M_{G_\gamma})}{\det(1-q\Fr_q\mid M_{G_\gamma})}\cdot\prod_{\alpha\in J_X\cup J_S-\{1,1\}}\det(1-\alpha\Fr_q\mid M_{G_\gamma}).
 \end{align*}
 
 \begin{exa}\label{exa:the sum of L-functions}
 Consider the case $X=\P^1_{\F_q}$ and $S=\{0,\infty\}$ with $\deg 0=\deg\infty=1$. Then one has $J_X=\emptyset$ and $J_S=\{1,1\}$, and so that
 \[
 L_S(M_{G_\gamma})=\frac{\det(1-\Fr_q\mid M_{G_\gamma})}{\det(1-q\Fr_q\mid M_{G_\gamma})}.
 \]
 Since we have assumed that $G$ is almost simple, it is shown in \cite[\S 9]{MR2843098} that
 \[
 \sum_{[\gamma]}\frac{\det(1-\Fr_q\mid M_{G_\gamma})}{\det(1-q\Fr_q\mid M_{G_\gamma})}=1.
 \]
 This is interesting because it states that the non-trivial sum of rational numbers equals 1.
 \end{exa}
 
 Our aim is to describe how the sum changes if the base field is changed from $\F_q$ to $\F_{q^m}$. However, this seems complicated because not only each $L_S(M_{G_\gamma})$ but also the set of semisimple classes which the sum is taken over will change under the base change. We start by stating a conjecture on the sum.
 
\begin{conj}\label{conj:the sum of L-functions}
Let $\mathcal{L}(G, m)$ denote the sum $\sum_{[\gamma]}L_S(M_{G_\gamma})$ for the base field $\F_{q^m}$ and $(X_m,S_m)$. Then the function $m\mapsto \mathcal{L}(G, m)$ is of Lefschetz type.
\end{conj}

\section{The sum of $L$-functions: $\SL_n$ case}\label{The sum of L-functions: SL_n case}
We first consider the case $G=\SL_n$. There is an explicit criterion of being conjugate under $\SL_n(\F_q)$.

\begin{lem}\label{lem:conjugacy in SL_n}
Two semisimple elements in $\SL_n(\F_q)$ are conjugate in $\SL_n(\F_q)$ if and only if they have the same characteristic polynomial over $\F_q$.
\end{lem}

\begin{prf}
It is clear that two semisimple elements in $\SL_n(\F_q)$ are conjugate in $\GL_n(\overline{\F}_q)$ if and only if they have the same characteristic polynomial over $\F_q$. However, the conjugacy in $\GL_n(\overline{\F}_q)$ implies the conjugacy in $\SL_n(\overline{\F}_q)$, so it suffices to show that the conjugacy in $\SL_n(\overline{\F}_q)$ implies the conjugacy in $\SL_n(\F_q)$. Since $\SL_n$ is simply connected and semisimple, we can apply \cite[II, 3.10]{MR268192}.
\end{prf}

By Lemma \ref{lem:conjugacy in SL_n}, we can classify the structures of centralizers of semisimple elements by seeing their characteristic polynomials. Let $\gamma\in G(\F_q)=\SL_n(\F_q)$ be a semisimple element. Then the characteristic polynomial of $\gamma$ is of the form
\[
P_1(x)^{a_1}\cdots P_r(x)^{a_r}
\]
where each $P_i(x)$ is an irreducible monic polynomial over $\F_q$ such that
\[
P_1(0)^{a_1}\cdots P_r(0)^{a_r}=(-1)^n.
\]
Put $d_i\coloneqq\deg P_i$ and call the weighted partition $n=\sum_{i=1}^r a_i d_i$ the \emph{type} of $\gamma$. We are going to show that the structure of the centralizer of a semisimple element is determined by its type.

Let $V\coloneqq \F_q^n$ be a vector space on which $G(\F_q)$ acts naturally. Then $V$ admits a decomposition
\[
V=V_{P_1}\oplus\cdots\oplus V_{P_r},\ \ \ V_{P_i}\coloneqq\Ker(P_i(\gamma)^{a_i}\colon V\to V)
\]
and we have
\begin{align*}
G_\gamma(\F_q)&=\GL(V)_{\gamma}\cap G(\F_q)\\
&=(\GL(V_{P_1})_{\gamma_1}\times\cdots\times\GL(V_{P_r})_{\gamma_r})\cap G(\F_q)
\end{align*}
where each $\gamma_i$ is the image of $\gamma$ in $\GL(V_{P_i})$.  For each $i$, the characteristic polynomial of $\gamma_i$ in $\GL(V_{P_i})$ is $P_i(x)^{a_i}$ and one has $\F_{q}[\gamma_i]\cong\F_{q^{d_i}}$. It follows that $\GL(V_{P_i})_{\gamma_i}\cong\GL_{a_i}(\F_{q^{d_i}})$ and that
\[
\GL_{n,\gamma}\cong \Res_{\F_{q^{d_1}}/\F_q}\GL_{a_1}\times\cdots\times\Res_{\F_{q^{d_r}}/\F_q}\GL_{a_r}.
\]
Hence, we obtain
\[
M_{\GL_{n,\gamma}}=\bigoplus_{i=1}^r(\Ind_{\F_{q^{d_i}}/\F_q}\Q\oplus\Ind_{\F_{q^{d_i}}/\F_q}\Q(-1)\oplus\cdots\oplus\Ind_{\F_{q^{d_i}}/\F_q}\Q(1-a_i))
\]
by Lemma \ref{lem:some properties of motives} and Example \ref{exa:motives}. On the other hand, the canonical isogeny $G\times \G_m\to\GL_n$ induces the isogeny
\[
G_{\gamma}\times\G_m\longrightarrow \GL_{n,\gamma}
\]
which gives the identity $M_{G_{\gamma}}\oplus\Q=M_{\GL_{n,\gamma}}$ by Lemma \ref{lem:some properties of motives}. It follows that
\begin{align*}
M_{G_\gamma}&=\left(\bigoplus_{i=1}^r\Ind_{\F_{q^{d_i}}/\F_q}\Q\right)/\Q\oplus\bigoplus_{i=1}^r(\Ind_{\F_{q^{d_i}}/\F_q}\Q(-1)\oplus\cdots\oplus\Ind_{\F_{q^{d_i}}/\F_q}\Q(1-a_i))\\
&=\left(\bigoplus_{i=1}^rW_{d_i}\right)/\Q\oplus\bigoplus_{i=1}^r(W_{d_i}(-1)\oplus\cdots\oplus W_{d_i}(1-a_i))
\end{align*}
where $W_{d_i}$ is the regular representation of $\Gal(\F_{q^{d_i}}/\F_q)$. Thus one can calculate as
\[
\det(1-t\Fr_q\mid M_{G_\gamma})=\frac{\prod_{i=1}^r(1-t^{d_i})}{1-t}\cdot\prod_{i=1}^r(1-(t q)^{d_i})\cdots(1-(t q^{a_i-1})^{d_i}).
\]
Consider the case $r=1$, namely, the case where the type of $\gamma$ is of the form $n=ad$. Then
\[
\det(1-t\Fr_q\mid M_{G_\gamma})=\frac{1-t^d}{1-t}\cdot(1-(t q)^d)\cdots(1-(t q^{a-1})^d)
\]
and especially we have
\[
\frac{\det(1-\Fr_q\mid M_{G_\gamma})}{\det(1-q\Fr_q\mid M_{G_\gamma})}=\frac{d\cdot (1-q^d)\cdots(1-q^{(a-1)d})}{\frac{1-q^d}{1-q}\cdot(1-q^{2d})\cdots(1-q^{ad})}=d\cdot\frac{1-q}{1-q^n}.
\]
On the other hand, if $r\geq 2$ then $\det(1-\Fr_q\mid M_{G_\gamma})=0$ since the factor
\[
\frac{\prod_{i=1}^r(1-t^{d_i})}{1-t}
\]
of $\det(1-t\Fr_q\mid M_{G_\gamma})$ is zero at $t=1$. In sum,

\begin{prop}\label{prop:the ratio of determinants}
We have
\[
\frac{\det(1-\Fr_q\mid M_{G_\gamma})}{\det(1-q\Fr_q\mid M_{G_\gamma})}=\begin{dcases}
d\cdot\frac{1-q}{1-q^n} & \text{if the type of $\gamma$ is $n=ad$}\\
0 & \text{otherwise}.
\end{dcases}
\]
In particular, $L_S(M_{G_\gamma})=0$ unless the type of $\gamma$ is of the form $n=ad$.
\end{prop}

In order to calculate the sum of $L$-functions, we need to count the number of possible characteristic polynomials of a given type. For each $d\mid n$, define $N_{n, d}(q)$ to be the number of irreducible monic polynomials $P(x)$ over $\F_q$ of degree $d$ such that $P(0)^{n/d}=(-1)^n$. If $\gamma$ is of type $n=ad$, denote the isomorphic class of $G_\gamma$ by $G_d$. Then $N_{n, d}(q)$ equals the number of semisimple classes corresponding to the type $n=ad$, and the sum of $L$-functions can be written as
\[
\sum_{[\gamma]}L_S(M_{G_\gamma})=\sum_{d\mid n}N_{n, d}(q)L_S(M_{G_d}).
\]
By Example \ref{exa:the sum of L-functions} and Proposition \ref{prop:the ratio of determinants}, we have the formula
\[
\sum_{d\mid n}dN_{n, d}(q)\frac{1-q}{1-q^n}=1.
\]

\begin{exa}\label{exa:N_{n, d}}
For any $n$, $N_{n,1}(q)$ is just the number of $n$-th roots of unity contained in $\F_q$: $N_{n,1}(q)=\#\mu_n(\F_q)$. If we assume that $n$ is a prime, then we have
\[
nN_{n,n}(q)=\frac{1-q^n}{1-q}-\#\mu_n(\F_q)
\]
by the above formula. In general, we expect that $N_{n, d}(q)$ is a polynomial in $q$ under some condition of the form $q\equiv i\bmod N$, but do not know the details.
\end{exa}

Let $\ell$ be a prime. The following is the answer to Conjecture \ref{conj:the sum of L-functions} for $G=\SL_\ell$.

\begin{thm}\label{thm:the sum of L-functions for SL_l}
The function $m\mapsto\mathcal{L}(\SL_\ell, m)$ is of Lefschetz type.
\end{thm}

\begin{prf}
The sum we want to calculate is 
\begin{align*}
\sum_{[\gamma]}L_S(M_{G_\gamma})&=N_{\ell,1}(q)L_S(M_{G_1})+N_{\ell,\ell}(q)L_S(M_{G_d})\\
&=N_{\ell,1}(q)\frac{1-q}{1-q^\ell}H_1(q)+\ell N_{\ell,\ell}(q)\frac{1-q}{1-q^\ell}H_\ell(q)
\end{align*}
where
\begin{align*}
H_1(q)&\coloneqq \prod_{\alpha\in J_X\cup J_S-\{1,1\}}\det(1-\alpha\Fr_q\mid M_{G_1})\\
&\,=\prod_{\alpha\in J_X\cup J_S-\{1,1\}}(1-\alpha q)(1-\alpha q^2)\cdots(1-\alpha q^{\ell-1}),\\
H_\ell(q)&\coloneqq \prod_{\alpha\in J_X\cup J_S-\{1,1\}}\det(1-\alpha\Fr_q\mid M_{G_\ell})\\
&\,=\prod_{\alpha\in J_X\cup J_S-\{1,1\}}(1+\alpha+\cdots+\alpha^{\ell-1}).
\end{align*}
which define polynomials $H_1(x),H_\ell(x)\in\Z[x]$. Then we have $H_1(\zeta_\ell)=H_\ell(\zeta_\ell)$ for $\zeta_\ell=e^{2\pi\sqrt{-1}/\ell}$ and thus $H_1(x)-H_\ell(x)$ can be written as $(1+x+\cdots+x^{\ell-1})R(x)$ for some $R(x)\in\Z[x]$. Therefore,
\begin{align*}
\sum_{[\gamma]}L_S(M_{G_\gamma})&=N_{\ell,1}(q)\frac{1-q}{1-q^\ell}(H_1(q)-H_\ell(q))+(N_{\ell,1}(q)+\ell N_{\ell,\ell}(q))\frac{1-q}{1-q^\ell}H_\ell(q)\\
&=\#\mu_\ell(\F_q)R(q)+H_\ell(q)\\
&=\begin{cases}
\ell R(q)+H_\ell(q) & \text{if } q\equiv 1\bmod\ell\\
R(q)+H_\ell(q) & \text{if } q\not\equiv 1\bmod\ell
\end{cases}
\end{align*}
by Example \ref{exa:N_{n, d}}. 

Consider the base change from $\F_q$ to $\F_{q^m}$. We have that $H_1(q)$, $H_\ell(q)$ and $R(q)$ are polynomials in the elements of $J\coloneqq\{q\}\cup J_X\cup J_S-\{1,1\}$ over $\Z$ by the above argument, so write as $H_1(J)=H_1(x)$, $H_\ell(J)=H_\ell(x)$ and $R(J)=R(x)$. Since $J$ is replaced by $J^m$ under the base change, we have
\[
\mathcal{L}(\SL_\ell, m)=\#\mu_\ell(\F_{q^m})R(J^m)+H_\ell(J^m).
\] 
It is clearly a Lefschetz type function in $m$.
\end{prf}

It is expected that Conjecture \ref{conj:the sum of L-functions} holds for $G=\SL_n$ in general. On this point, we have the following result:

\begin{thm}\label{thm:the sum of L-functions for SL_n}
Let $J$ be a finite set of indeterminates corresponding to the multiset $J_X\cup J_S-\{1,1\}$. Then there exists a polynomial $\Pcal(x,J)\in\Z[x,J]$ such that 
\[
\mathcal{L}(\SL_n,m)=\Pcal(q^m,J_X^m\cup  J_S^m-\{1,1\})
\]
for any $m\in\Z_{>0}$ such that $\gcd(n,q^m-1)=1$. 
\end{thm}

We need some preparations:

\begin{lem}\label{lem:N_{n, d}(q)}
Assume that $\gcd(n, q-1)=1$. Then for any $d\mid n$, we have
\[
d N_{n, d}(q)=d N_{d, d}(q)=\sum_{e\mid d}\mu(e)\frac{q^{d/e}-1}{q-1}.
\]
Here, $\mu$ is the Möbius function.
\end{lem}

\begin{prf}
Let $P(x)$ be an irreducible monic polynomial over $\F_q$ of degree $d$. Then the condition $P(0)^{n/d}=(-1)^n$ is equivalent to the condition $P(0)=(-1)^d$ due to our assumption which implies $\mu_{n/d}(\F_q)=\{1\}$. This shows the first equality. Now we have the equality
\[
\sum_{d\mid n}d N_{d, d}(q)=\frac{q^n-1}{q-1}
\]
for any $d\mid n$, by the formula just before Example \ref{exa:N_{n, d}}. Then the Möbius inversion formula gives the second equality.
\end{prf}

\begin{lem}\label{lem:H_{n, d}(x,J)}
For $n, d\in\Z_{>0}$ such that $d\mid n$, define the polynomial $H_{n, d}(x,J)\in\Z[x,J]$ as 
\[
H_{n, d}(x,J)\coloneqq \prod_{\alpha\in J}\frac{1}{1-\alpha}(1-\alpha^d)(1-(\alpha x)^d)(1-(\alpha x^2)^d)\cdots(1-(\alpha x^{n/d-1})^d).
\]
Here, $J$ is a finite set of indeterminates. 
\begin{enumerate}
\item For any $m\in\Z_{>0}$, 
\[
H_{mn, md}(x,J)=H_{n, d}(x^m,J^m)\prod_{\alpha\in J}\frac{1-\alpha^m}{1-\alpha}.
\]
\item For any $c\mid n$,
\[
H_{n, d}(\zeta_c,J)=\prod_{\alpha\in J}\frac{(1-\alpha^{\lcm(c, d)})^{n/\lcm(c, d)}}{1-\alpha}.
\]
\end{enumerate}
\end{lem}

\begin{prf}
One can check (1) easily from the definition. For (2), put $r\coloneqq c/{\gcd(c, d)}$ so that $\zeta_c^d\in\C^\times$ is a primitive $r$-th root of unity. It follows that
\begin{align*}
H_{n, d}(\zeta_c,J)&=\prod_{\alpha\in J}\frac{1}{1-\alpha}(1-\alpha^d)(1-\zeta_c^d\alpha^d)(1-(\zeta_c^d)^2\alpha^d)\cdots(1-(\zeta_c^d)^{n/d-1}\alpha^d)\\
&=\prod_{\alpha\in J}\frac{1}{1-\alpha}(1-\alpha^{rd})^{n/(rd)}\\
&=\prod_{\alpha\in J}\frac{(1-\alpha^{\lcm(c, d)})^{n/\lcm(c, d)}}{1-\alpha}.
\end{align*}
\end{prf}

Now we are going to prove Theorem \ref{thm:the sum of L-functions for SL_n}. We have
\[
\mathcal{L}(\SL_n,m)=\sum_{d\mid n}d N_{n, d}(q^m)\frac{q^m-1}{(q^m)^n-1}H_{n, d}(q^m,J_X^m\cup J_S^m-\{1,1\}).
\]
by Proposition \ref{prop:the ratio of determinants}. We put
\[
M_{n, d}(x)\coloneqq \sum_{e\mid d}\mu(e)\frac{x^{d/e}-1}{x^n-1}
\]
for any $n, d\in\Z_{>0}$ such that $d\mid n$. Then one has
\[
M_{n, d}(q^m)=d N_{n, d}(q^m)\frac{q^m-1}{(q^m)^n-1}
\]
by Lemma \ref{lem:N_{n, d}(q)} under the assumption $\gcd(n, q^m-1)=1$. Therefore, it suffices to show that the rational function
\[
\Pcal(x,J)\coloneqq\sum_{d\mid n}H_{n, d}(x,J)M_{n, d}(x)
\]
lies in $\Z[x,J]$. In fact, we have the following stronger result:

\begin{prop}
For any $n',d'\in\Z_{>0}$ such that $d'\mid n'$ and $d'$ is prime to $n$,
\[
\sum_{d\mid n}H_{n'n, d'd}(x,J)M_{n, d}(x)
\]
lies in $\Z[x,J]$.
\end{prop}

\begin{prf}
If $n=1$, then the assertion is trivial. Assume that the proposition is true for $n=N$, and we are going to prove it for $n=\ell^kN$, where $\ell$ is a prime not dividing $N$, and $k$ is a positive integer. Fix $n',d'\in\Z_{>0}$ such that $d'\mid n'$ and $d'$ is prime to $\ell^k N$. We have
\begin{align*}
\sum_{d\mid \ell^kN}H_{n'\ell^kN, d'd}(x,J)M_{\ell^kN, d}(x)=\sum_{i=0}^k\sum_{d\mid N}H_{n'\ell^kN,d'\ell^id}(x,J)M_{\ell^kN,\ell^id}(x)
\end{align*}
and
\begin{align*}
M_{\ell^kN,\ell^id}(x)&=\sum_{e\mid \ell^id}\mu(e)\frac{x^{\ell^id/e}-1}{x^{\ell^kN}-1}=\sum_{j=0}^i\sum_{e\mid d}\mu(\ell^j e)\frac{x^{\ell^id/(\ell^je)}-1}{x^{\ell^kN}-1}\\
&=\sum_{j=0}^i\mu(\ell^j)\sum_{e\mid d}\mu(e)\frac{(x^{\ell^{i-j}})^{d/e}-1}{(x^{\ell^{i-j}})^N-1}\cdot\frac{x^{\ell^{i-j}N}-1}{x^{\ell^kN}-1}\\
&=\begin{dcases}
\frac{x^N-1}{x^{\ell^kN}-1}M_{N, d}(x) & \text{if } i=0 \\
\frac{x^{\ell^iN}-1}{x^{\ell^kN}-1}M_{N, d}(x^{\ell^i})-\frac{x^{\ell^{i-1}N}-1}{x^{\ell^kN}-1}M_{N, d}(x^{\ell^{i-1}}) & \text{if } i\geq 1.
\end{dcases}
\end{align*}
Hence,
\begin{align*}
&\sum_{d\mid \ell^kN}H_{n'\ell^kN, d'd}(x,J)M_{\ell^kN, d}(x)\\
&\quad=\sum_{d\mid N}H_{n'\ell^kN, d'd}(x,J)\frac{x^N-1}{x^{\ell^kN}-1}M_{N, d}(x)\\
&\quad\quad+\sum_{i=1}^k\sum_{d\mid N}H_{n'\ell^kN,d'\ell^id}(x,J)\left(\frac{x^{\ell^iN}-1}{x^{\ell^kN}-1}M_{N, d}(x^{\ell^i})-\frac{x^{\ell^{i-1}N}-1}{x^{\ell^kN}-1}M_{N, d}(x^{\ell^{i-1}})\right)\\
&\quad=\sum_{i=0}^{k}\frac{x^{\ell^iN}-1}{x^{\ell^kN}-1}\sum_{d\mid N}H_{n'\ell^kN,d'\ell^id}(x,J)M_{N, d}(x^{\ell^i})\\
&\quad\quad-\sum_{i=0}^{k-1}\frac{x^{\ell^iN}-1}{x^{\ell^kN}-1}\sum_{d\mid N}H_{n'\ell^kN,d'\ell^{i+1}d}(x,J)M_{N, d}(x^{\ell^i}).
\end{align*}
Here, if we put
\begin{align*}
P_i(x,J)&\coloneqq\sum_{d\mid N}H_{n'\ell^kN,d'\ell^id}(x,J)M_{N, d}(x^{\ell^i})\text{ for } i=0,\ldots, k,\\
Q_i(x,J)&\coloneqq\sum_{d\mid N}H_{n'\ell^kN,d'\ell^{i+1}d}(x,J)M_{N, d}(x^{\ell^i})\text{ for } i=0,\ldots,k-1
\end{align*}
then one sees that
\begin{align*}
&P_i(x,J)=\sum_{d\mid N}H_{n'\ell^{k-i}N, d'd}(x^{\ell^i},J^{\ell^i})M_{N, d}(x^{\ell^i})\cdot\prod_{\alpha\in J}\frac{1-\alpha^{\ell^i}}{1-\alpha},\\
&Q_i(x,J)=\sum_{d\mid N}H_{n'\ell^{k-i}N,d'\ell d}(x^{\ell^i},J^{\ell^i})M_{N, d}(x^{\ell^i})\cdot\prod_{\alpha\in J}\frac{1-\alpha^{\ell^i}}{1-\alpha}
\end{align*}
due to Lemma \ref{lem:H_{n, d}(x,J)} (1), and they lie in $\Z[x,J]$ by the induction hypothesis. Now we have
\begin{align*}
&\sum_{d\mid \ell^kN}H_{n'\ell^kN, d'd}(x,J)M_{\ell^kN, d}(x)\\
&\quad=\sum_{i=0}^{k}\frac{x^{\ell^iN}-1}{x^{\ell^kN}-1}P_i(x,J)-\sum_{i=0}^{k-1}\frac{x^{\ell^iN}-1}{x^{\ell^kN}-1}Q_i(x,J)\\
&\quad=\sum_{i=0}^{k-1}\frac{x^{\ell^iN}-1}{x^{\ell^kN}-1}(P_i(x,J)-Q_i(x,J))+P_k(x,J).
\end{align*}
Let $\Phi_m(x)$ be the cyclotomic polynomial for $m\in\Z_{>0}$. Then one has
\[
\frac{x^{\ell^iN}-1}{x^{\ell^kN}-1}=\frac{1}{\prod_{c\mid\ell^kN, c\nmid\ell^iN}\Phi_c(x)}=\frac{1}{\prod_{\ell^{i+1}\mid c\mid \ell^k N}\Phi_c(x)}.
\]
Therefore, it suffices to show that $P_i(\zeta_c,J)-Q_i(\zeta_c,J)=0$ for $i=0,\ldots, k-1$ and $c$ such that $\ell^{i+1}\mid c\mid \ell^k N$. We have
\[
M_{N, d}(x^{\ell^i})=\sum_{e\mid d}\mu(e)\frac{x^{\ell^id/e}-1}{x^{\ell^iN}-1}
\]
and thus $M_{N, d}(\zeta_c^{\ell^i})<\infty$ if $\ell^{i+1}\mid c$. So it suffices to show that $H_{n'\ell^kN,d'\ell^id}(\zeta_c,J)=H_{n'\ell^kN,d'\ell^{i+1}d}(\zeta_c,J)$. However, this is an immediate consequence of Lemma \ref{lem:H_{n, d}(x,J)} (2) since one has $\lcm(c, d'\ell^id)=\lcm(c, d'\ell^{i+1}d)$ for such $i$ and $c$.
\end{prf}

\section{The sum of $L$-functions: $\Sp_{2n}$ case}\label{The sum of L-functions: Sp_2n case}
Next we work on the case $G=\Sp_{2n}$. There is a similar criterion to Lemma \ref{lem:conjugacy in SL_n}

\begin{lem}
Two semisimple elements in $\Sp_{2n}(\F_q)$ are conjugate in $\Sp_{2n}(\F_q)$ if and only if they have the same characteristic polynomial over $\F_q$. 
\end{lem}

\begin{prf}
One has the fact that two elements in $\Sp_{2n}(\overline{\F}_q)$ are conjugate in $\Sp_{2n}(\overline{\F}_q)$ if and only if they are in $\GL_{2n}(\overline{\F}_q)$ \cite[IV, 2.15]{MR268192}, so it suffices to show that the conjugacy in $\Sp_{2n}(\overline{\F}_q)$ implies the conjugacy in $\Sp_{2n}(\F_q)$. Since $\Sp_{2n}$ is simply connected and semisimple, we can apply \cite[II, 3.10]{MR268192}.
\end{prf}

Let $\gamma\in G(\F_q)=\Sp_{2n}(\F_q)$ be a semisimple element. Since the multiset of eigenvalues of $\gamma$ is stable under taking inverse, the characteristic polynomial of $\gamma$ is of the form
\[
     (x-1)^{2a_{+}}(x+1)^{2a_{-}}P_1(x)^{b_1}\cdots P_r(x)^{b_r}(Q_1(x)Q_1^*(x))^{c_1}\cdots(Q_s(x)Q_s^*(x))^{c_s}
\]
where
\begin{itemize}
    \item[\textbullet] $x\pm 1$, $P_i(x)$, $Q_j(x)$ are distinct irreducible monic polynomials over $\F_q$.
    \item[\textbullet] $a_+,a_-\in \Z_{\geq 0}$, $b_i, c_j\in\Z_{>0}$.
    \item[\textbullet] $P_i(x)$ are self-reciprocal ($P_i=P_i^*$).
    \item[\textbullet] $Q_j(x)$ are not self-reciprocal ($Q_j\neq Q_j^*$). 
\end{itemize}
Here, for any polynomial $R(x)$ with non-zero constant term, the reciprocal polynomial $R^*(x)$ of $R(x)$ is defined as $R^*(x)\coloneqq R(0)^{-1}x^{\deg R}R(x^{-1})$. Note that each $P_i(x)$ is of even degree. Indeed, there are only two self-reciprocal irreducible monic polynomial of odd degree: $x\pm 1$. If $q$ is even, then we suppose $a_-$ to be zero.

Let $V=\F_q^{2n}$ be a vector space with a symplectic form $\langle\,,\,\rangle\colon V\times V\to \F_q$ so that $G(\F_q)=\mathrm{Sp}(V)$ is defined by $\langle\,,\,\rangle$. Then $V$ admits a decomposition
\[
V=V_{x-1}\oplus V_{x-1}\oplus\bigoplus_{i=1}^rV_{P_i}\oplus\bigoplus_{j=1}^s(V_{Q_j}\oplus V_{Q_j^*})
\]
where $V_{x\pm 1}\coloneqq\Ker((\gamma\pm1)^{2a_\pm}\colon V\to V)$,  $V_{P_i}\coloneqq\Ker(P_i(\gamma)^{b_i}\colon V\to V)$, $V_{Q_j}\coloneqq\Ker(Q_j(\gamma)^{c_j}\colon V\to V)$. It is obvious that
\begin{align*}
G_\gamma(\F_q)\subset\GL(V_{x-1})\times\GL(V_{x+1})\times\prod_{i=1}^r\GL(V_{P_i})\times\prod_{j=1}^s(\GL(V_{Q_j})\times\GL(V_{Q_j^*})).
\end{align*}

\begin{lem}\label{lem:restrictions of symplectic forms}
\begin{enumerate}
\item The pairing $\langle\,,\,\rangle\colon V\times V\to \F_q$ induces symplectic forms on $V_{x\pm 1}$, $V_{P_i}$, and $V_{Q_j}\oplus V_{Q_j^*}$.
\item The subspaces $V_{Q_j}$ and $V_{Q_j^*}$ are isotropic with respect to $\langle\,,\,\rangle$.
\end{enumerate}
\end{lem}

\begin{prf}
For (1), it suffices to show that the restrictions are non-degenerate. Let $R(x)$ be any of $(x\pm 1)^{a_\pm}$, $P_i(x)^{b_i}$, $(Q_j(x)Q_j^*(x))^{c_j}$ and put $V_R\coloneqq \Ker(R(\gamma)\colon V\to V)$. We will prove that  $V_R\cap V_R^\perp=0$. Write the above decomposition as $V=V_R\oplus W$. Then $R(\gamma)\colon V\to V$ induces a bijective map $R(\gamma)|_W\colon W\to W$ and thus we have
\begin{align*}
\langle V_R,W\rangle=\langle V_R,R(\gamma)W\rangle=\langle R(\gamma^{-1})V_R,W\rangle=\langle R(0)\gamma^{-\deg R}R(\gamma)V_R,W\rangle=0.
\end{align*}
Here we used $\langle-,\gamma(-)\rangle=\langle\gamma^{-1}(-),-\rangle$ and $R=R^*$. This shows
\[
V_R\cap V_R^\perp\subset W^\perp\cap V_R^\perp=(W+V_R)^\perp=0.
\]

For (2), let $R(x)$ be either $Q_j(x)^{c_j}$ or $Q_j^*(x)^{c_j}$. Then 
\[
\langle V_R,V_R \rangle=\langle V_R,R^*(\gamma)V_R\rangle=\langle R(0)\gamma^{-\deg R}R(\gamma)V_R,V_R\rangle=0
\]
by the same argument as above, and the assertion follows.
\end{prf}

By Lemma \ref{lem:restrictions of symplectic forms} (1), we can reduce to the cases where the characteristic polynomial of $\gamma$ is of the form
\begin{enumerate}
    \item[(a)] $(x\pm 1)^{2a_{\pm}}$,
    \item[(b)] $P(x)^b$ with $P=P^*$,
    \item[(c)] $(Q(x)Q^*(x))^c$ with $Q\neq Q^*$.
\end{enumerate}
The case (a) is obvious: $G_\gamma=\mathrm{Sp}_{2n}$. Now consider the case (b). We put $d\coloneqq\deg P/2$. Then $V$ can be regarded as a vector space of dimension $b$ over $\F_q[\gamma]\cong\F_{q^{2d}}$ and $G_\gamma(\F_q)$ corresponds to the set of $g\in\mathrm{Aut}_{\F_{q^{2d}}}(V)$ preserving the symplectic form $\langle\,,\,\rangle$. Let $\sigma\in\Gal(\F_{q^{2d}}/\F_q)$ be the involution defined by $\sigma(\gamma)=\gamma^{-1}$. Then $\langle\,,\,\rangle$ is a $\sigma$-skew pairing, namely, it satisfies $\langle g x, y\rangle=\langle x,\sigma(g)y\rangle$ for $g\in \F_{q^{2d}}$, $x, y\in V$. Define $\llparenthesis\,,\,\rrparenthesis\colon V\times V\to \F_{q^{2d}}$ as follows. Let $\{u_i\}_{1\leq i\leq 2d}$ be an $\F_q$-basis of $\F_{q^{2d}}$ and $\{v_i\}_{1\leq i\leq 2d}$ be its dual basis with respect to $\Tr_{\F_{q^{2d}}/\F_q}$. Define
\[
\llparenthesis x, y\rrparenthesis\coloneqq\sum_{i=1}^{2d}v_i\langle x, u_i y\rangle\ \ \ (x, y\in V).
\]
Then it is shown in \cite[Lemma 1.1.4.5]{5c715bee-b469-388d-a7b5-82794e13039c} that $\llparenthesis\,,\,\rrparenthesis$ is a skew-Hermitian form over $\F_{q^{2d}}$ such that $\Tr_{\F_{q^{2d}}/\F_q}\circ \llparenthesis\,,\,\rrparenthesis=\langle\,,\,\rangle$, and preserving $\langle\,,\,\rangle$ is equivalent to preserving $\llparenthesis\,,\,\rrparenthesis$. So we have
\[
G_\gamma(\F_q)=\{g\in\GL_b(\F_{q^{2d}})\mid \llparenthesis g x, g y\rrparenthesis=\llparenthesis x, y\rrparenthesis,\ \forall x, y\in V\}\cong \U_b(\F_{q^d})
\]
where $\U_b$ is the unitary group of rank $b$ with respect to the quadratic extension $\F_{q^{2d}}/\F_{q^d}$. It follows that $G_\gamma\cong\Res_{\F_{q^{d}}/\F_q}\U_b$. 

Next we consider the case (c). Put $e\coloneqq\deg Q$. We have 
\[
G_\gamma(\F_q)\subset (\GL(V_Q)\times\GL(V_{Q^*}))\cap G(\F_q)\cong \GL(V_Q).
\]
Indeed, $V_Q$ is a maximal isotropic subspace of $V=V_Q\oplus V_{Q^*}$ by Lemma \ref{lem:restrictions of symplectic forms} (2), and so we may assume that $\langle\,,\,\rangle$ is defined by $J=\left(\begin{smallmatrix}0&I_n\\-I_n&0\end{smallmatrix}\right)$ by choosing bases of $V_Q$ and $V_{Q^*}$. Then $\,^t\!\left(\begin{smallmatrix}A&0\\0&\,B\end{smallmatrix}\right)J\left(\begin{smallmatrix}A&0\\0&B\end{smallmatrix}\right)=J$ implies $B=\,^t\!A^{-1}$ for $A\in\GL(V_Q)$, $B\in\GL(V_{Q^*})$. Let $\overline{\gamma}$ denote the image of $\gamma$ in $\mathrm{GL}(V_Q)$. Then $V_Q$ becomes a vector space of dimension $c$ over $\F_q[\overline{\gamma}]\cong\F_{q^e}$. By a similar argument as in (a), 
\[
G_\gamma(\F_q)\cong \GL(V_Q)_{\overline{\gamma}}\cong\mathrm{Aut}_{\F_{q^e}}(V)\cong \GL_c(\F_{q^e}). 
\]
It follows that $G_\gamma\cong \Res_{\F_{q^e}/\F_q}\GL_c$.

In sum, we have the general result:
\begin{prop}\label{prop:centralizers in Sp2n}
    The centralizer $G_\gamma$ is of the form 
    \[
    \Sp_{2a_+}\times \Sp_{2a_-}\times\prod_{i=1}^r\Res_{\F_{q^{d_i}}/\F_q}\U_{b_i}\times \prod_{j=1}^s\Res_{\F_{q^{e_j}}/\F_q}\GL_{c_j}
    \]
   where $a_\pm\in\Z_{\geq 0}$, $b_i, c_j, d_i, e_j\in\Z_{>0}$. 
\end{prop}

If $G_\gamma$ is of the form in Proposition \ref{prop:centralizers in Sp2n}, we call the tuple 
\[
(a_+,a_-;(d_i, b_i)_{1\leq i\leq r};(e_j, c_j)_{1\leq j\leq s})
\]
the \emph{type} of $\gamma$. We consider two types $\tau=(a_+,a_-;((d_i, b_i))_{1\leq i\leq r};((e_j, c_j))_{1\leq j\leq s})$ and $\tau'=(a'_+,a'_-;((d'_i, b'_i))_{1\leq i\leq r'};((e'_j, c'_j))_{1\leq j\leq s'})$ to be the same if $r=r'$, $s=s'$, and $(a_+,a_-)$, $((d_i, b_i))_{1\leq i\leq r}$, $((e_j, c_j))_{1\leq j\leq s}$ correspond to $(a'_+,a'_-)$, $((d'_i, b'_i))_{1\leq i\leq r'}$, $((e'_j, c'_j))_{1\leq j\leq s'}$ up to ordering, respectively. Then Proposition \ref{prop:centralizers in Sp2n} implies that the isomorphic class of $G_\gamma$ depends only on the type of $\gamma$. 

Next we determine the motive $M_{G_\gamma}$ of $G_\gamma$ and its $L$-function. We have
\begin{align*}
&M_{\Sp_{2a}}=\Q(-1)\oplus\Q(-3)\oplus\cdots\oplus\Q(1-2a),\\
&M_{\Res_{\F_{q^d}/\F_q}\U_b}=W_d[\sigma]\oplus W_d(-1)\oplus W_d[\sigma](-2)\oplus\cdots\oplus W_d[\sigma^b](1-b),\\
&M_{\Res_{\F_{q^e}/\F_q}\GL_c}=W_e\oplus W_e(-1)\oplus\cdots\oplus W_e(1-c)
\end{align*}
by Lemma \ref{lem:some properties of motives} and Example \ref{exa:motives}. Here, $W_d=\Ind_{\F_{q^d}/\F_q}\Q$ (resp.\ $W_e=\Ind_{\F_{q^e}/\F_q}\Q$) denote the regular representation of $\Gal(\F_{q^d}/\F_q)$ (resp.\ $\Gal(\F_{q^e}/\F_q)$), 
and $W_e[\sigma]=\Ind_{\F_{q^e}/\F_q}\Q[\sigma]$ is the representation on which $\Fr_q$ acts as $\left(\begin{smallmatrix}&1&&\\&&\ddots&\\&&&1\\-1&&&\end{smallmatrix}\right)$. Thus one can calculate as
\begin{align*}
    &\det(1-t\Fr_q\mid M_{\mathrm{Sp}_{2a}})=(1-t q)(1-t q^3)\cdots(1-t q^{2a}),\\
    &\det(1-t\Fr_q\mid M_{\mathrm{Res}_{\F_{q^d}/\F_q}\mathrm{U}_b})=(1+t^d)(1-(t q)^d)(1+(t q^2)^d)\cdots(1-(-1)^b(t q^{b-1})^d),\\
    &\det(1-t\Fr_q\mid M_{\mathrm{Res}_{\F_{q^e}/\F_q}\mathrm{GL}_c})=(1-t^e)(1-(t q)^e)\cdots(1-(t q^{c-1})^e).
\end{align*}
In general, $\det(1-t\Fr_q\mid M_{G_\gamma})$ is a finite product of polynomials of the above forms. Note that the third is zero at $t=1$, which implies that $L_S(M_{\Res_{\F_{q^e}/\F_q}\GL_c})$ vanishes. It follows that $L_S(M_{G_\gamma})$ is zero unless the characteristic polynomial of $\gamma$ is of the form 
\[
(x-1)^{2a_{+}}(x+1)^{2a_{-}}P_1(x)^{a_1}\cdots P_r(x)^{a_r}.
\]
Therefore, it is enough to deal with only such $\gamma$'s as long as we consider the sum of $L$-functions. Let $N_\tau(q)$ be the number of semisimple classes of type $\tau$. If $\gamma$ is of type $\tau$, we denote the isomorphic class of $G_\gamma$ by $G_\tau$. Then the sum of $L$-functions can be written as
\[
\sum_{[\gamma]}L_S(M_{G_\gamma})=\sum_{\tau}N_\tau(q)L_S(M_{G_\tau}).
\]
To describe $N_\tau(q)$, we need the following lemma:
\begin{lem}
The number of self-reciprocal irreducible monic polynomials over $\F_q$ of degree $2n$ is
\[
S_{2n}(q)=\begin{dcases}
\frac{1}{2n}(q^n-1) & \text{if } q\text{ odd and } n=2^s \\
\frac{1}{2n}\sum_{2\nmid d\mid n}\mu(d)q^{n/d} & \text{otherwise}.
\end{dcases}
\]
\end{lem}

\begin{prf}
See \cite[Theorem 2]{MR215815}.
\end{prf}

Assume that $\tau=(a_+,a_-;((d_i, b_i))_{1\leq i\leq r};)$ and define $k_{(d, b)}\coloneqq \#\{i\mid (d_i, b_i)=(d, b)\}$, $k_d\coloneqq \sum_b k_{(d, b)}=\#\{i\mid d_i=d\}$ for each $(d, b)\in\Z_{>0}\times\Z_{>0}$. Then we have
\[
N_\tau(q)=\begin{dcases}
\prod_d \frac{S_{2d}(q)(S_{2d}(q)-1)\cdots(S_{2d}(q)-k_d+1)}{\prod_b k_{(b, d)}!} & \text{if }a_+=a_- \text{ or } q \text{ even}\\
2\prod_d \frac{S_{2d}(q)(S_{2d}(q)-1)\cdots(S_{2d}(q)-k_d+1)}{\prod_b k_{(b, d)}!} & \text{if }a_+\neq a_- \text{ and } q \text{ odd}.
\end{dcases}
\]
Indeed, the set of semisimple classes of a given type corresponds to the set of characteristic polynomials of that type, except that interchanging $a_+$ and $a_-$ gives the same type but a different polynomial if $a_+\neq a_-$ and $q$ is odd. In particular, $N_\tau(q)$ is a polynomial in $q$ over $\Q$. 

We are going to observe the sum of $L$-functions in some $n$ small cases. We are done for $n=1$ because $\mathrm{Sp}_2=\mathrm{SL}_2$. Consider the case $G=\mathrm{Sp}_4$. Then there are six types of semisimple classes if $q$ is odd, and five types if $q$ is even. (See Table 1, 2 for details.)

\renewcommand{\arraystretch}{1.3}
\begin{table}[ht]
\centering
\caption{types of semisimple classes and related data for $\Sp_4$, $q$ odd }

\begin{tabular}{|c|c|c|c|}\hhline{|-|-|-|-|}
    type $\tau$ & $G_\tau$ & $N_\tau(q)$ & $\det(1-t\Fr_q\mid M_{G_\tau})$ \\ \hhline{|=|=|=|=|}
    $\tau_1=(2,0;;)$ & $\Sp_4$ & 2 & $(1-t q)(1-t q^3)$\\ \hhline{|-|-|-|-|}
    $\tau_2=(1,1;;)$ & $\Sp_2\times\Sp_2$ & $1$ & $(1-t q)^2$\\ \hhline{|-|-|-|-|}
    $\tau_3=(1,0;(1,1);)$ & $\Sp_2\times \U_1$ & $q-1$ & $(1-t q)(1+t)$\\ \hhline{|-|-|-|-|}
    $\tau_4=(;(1,2);)$ & $\U_2$ & $(q-1)/2$ & $(1+t)(1-t q)$ \\ \hhline{|-|-|-|-|}
    $\tau_5=(;(1,1),(1,1);)$ & $\U_1\times\U_1$ & $(q-1)(q-3)/8$ & $(1+t)^2$ \\ \hhline{|-|-|-|-|}
    $\tau_6=(;(2,1);)$ & $\Res_{\F_{q^2}/\F_q}\U_1$ & $(q^2-1)/4$ & $1+t^2$ \\ \hhline{|-|-|-|-|}
\end{tabular}

\caption{types of semisimple classes and related data for $\Sp_4$, $q$ even }
\begin{tabular}{|c|c|c|c|}\hhline{|-|-|-|-|}
    type $\tau'$ & $G_{\tau'}$ & $N_{\tau'}(q)$ & $\det(1-t\Fr_q\mid M_{G_{\tau'}})$ \\ \hhline{|=|=|=|=|}
    $\tau'_1=(2;;)$ & $\Sp_4$ & 1 & $(1-t q)(1-t q^3)$\\ \hhline{|-|-|-|-|}
    $\tau'_3=(1;(1,1);)$ & $\Sp_2\times \U_1$ & $q/2$ & $(1-t q)(1+t)$\\ \hhline{|-|-|-|-|}
    $\tau'_4=(;(1,2);)$ & $\U_2$ & $q/2$ & $(1+t)(1-t q)$ \\ \hhline{|-|-|-|-|}
    $\tau'_5=(;(1,1),(1,1);)$ & $\U_1\times\U_1$ & $q(q-2)/8$ & $(1+t)^2$ \\ \hhline{|-|-|-|-|}
    $\tau'_6=(;(2,1);)$ & $\Res_{\F_{q^2}/\F_q}\U_1$ & $q^2/4$ & $1+t^2$ \\ \hhline{|-|-|-|-|}
\end{tabular}
\end{table}
\renewcommand{\arraystretch}{1}

\begin{thm}\label{thm:the sum of L-functions for Sp_4}
   The function $m\mapsto\mathcal{L}(\mathrm{Sp}_4,m)$ is of Lefschetz type. 
\end{thm}

\begin{prf}
First, we assume that $q$ is odd. We define the rational functions $R_{\tau_i}(x)\in\Z(x)$ as
\begin{align*}
&R_{\tau_1}(x)=\frac{2(1+x+x^2)}{(1+x)^2(1+x^2)},\ \ R_{\tau_2}(x)=\frac{1}{(1+x)^2},\ \ R_{\tau_3}(x)=\frac{2(x-1)}{(1+x)^2},\\
&R_{\tau_4}(x)=\frac{x-1}{(1+x)^2},\ \ R_{\tau_5}(x)=\frac{(x-1)(x-3)}{2(1+x)^2},\ \ R_{\tau_6}(x)=\frac{x^2-1}{2(1+x^2)}
\end{align*}
and the polynomials $H_{\tau_i}(x)\in\Z[x]$ as
\begin{align*}
&H_{\tau_1}(x)=\prod_{\alpha\in J}(1-\alpha x)(1-\alpha x^3),\ \ H_{\tau_2}(x)=\prod_{\alpha\in J}(1-\alpha x)^2,\\
&H_{\tau_3}(x)=\prod_{\alpha\in J}(1+\alpha)(1-\alpha x),\ \ H_{\tau_4}(x)=\prod_{\alpha\in J}(1+\alpha)(1-\alpha x),\\
&H_{\tau_5}(x)=\prod_{\alpha\in J}(1+\alpha)^2,\ \ H_{\tau_6}(x)=\prod_{\alpha\in J}(1+\alpha^2).
\end{align*}
Here we put $J\coloneqq J_X\cup J_S-\{1,1\}$. Then one can see that
\begin{align*}
&R_{\tau_i}(q)=N_{\tau_i}(q)\frac{\det(1-\Fr_q\mid M_{G_{\tau_i}})}{\det(1-q\Fr_q\mid M_{G_{\tau_i}})},\\
&H_{\tau_i}(q)=\prod_{\alpha\in J}\det(1-\alpha \Fr_q\mid M_{G_{\tau_i}})
\end{align*}
for $i=1,\ldots, 6$ using Table 1. We put 
\[
\Pcal(x)\coloneqq \sum_{i=1}^6R_{\tau_i}(x)H_{\tau_i}(x)\in \Z(x)
\]
so that
\[
\Pcal(q)=\sum_{[\gamma]}L_S(M_{G_\gamma})=\sum_{i=1}^6N_{\tau_i}(q)L_S(M_{G_{\tau_i}}).
\]
Now we are going to prove that $\Pcal(x)\in\Z[x]$. One can see that $\Pcal(x)$ lies in $(2(1+x)^2(1+x^2))^{-1}\Z[x]$. We put $\Pcal_1(x)\coloneqq 2\Pcal(x)$, $\Pcal_2(x)\coloneqq (1+x)^2\Pcal(x)$, $\Pcal_3(x)\coloneqq (1+x^2)\Pcal(x)$ and show that
\begin{enumerate}
\item[(a)] $\Pcal_1(x)\equiv 0\bmod{2}$ in $((1+x)^2(1+x^2))^{-1}\Z[x]$,
\item[(b)] $\Pcal_2(x)\equiv 0\bmod{(1+x)^2}$ in $(2(1+x^2))^{-1}\Z[x]$,
\item[(c)] $\Pcal_3(x)\equiv 0\bmod{1+x^2}$ in $(2(1+x)^2)^{-1}\Z[x]$.
\end{enumerate}

For (a), one has
\begin{align*}
\Pcal_1(x)&\equiv 2(R_{\tau_5}(x)H_{\tau_5}(x)+R_{\tau_6}(x)H_{\tau_6}(x))\\
&=\frac{(x-1)(x-3)}{(1+x)^2}\prod_\alpha(1+\alpha)^2+\frac{x^2-1}{1+x^2}\prod_\alpha(1+\alpha^2)\\
&\equiv \prod_\alpha(1+\alpha)^2+\prod_\alpha(1+\alpha^2)\bmod{2}.
\end{align*}
It is zero modulo $2\overline{\Z}$, and thus is zero modulo $2\Z$. 

For (c), it suffices to show that $\Pcal_3(\sqrt{-1})=0$. One has
\begin{align*}
\Pcal_3(x)&\equiv(1+x^2)(R_{\tau_1}(x)H_{\tau_1}(x)+R_{\tau_6}(x)H_{\tau_6}(x))\\
&=\frac{2(1+x+x^2)}{(1+x)^2}\prod_\alpha(1-\alpha x)(1-\alpha x^3)+\frac{x^2-1}{2}\prod_\alpha(1+\alpha^2)\bmod{1+x^2}
\end{align*}
and so that 
\begin{align*}
\Pcal_3(\sqrt{-1})=1\cdot\prod_\alpha(1+\alpha^2)+(-1)\cdot\prod_\alpha(1+\alpha^2)=0.
\end{align*}

For (b), it suffices to show that $\Pcal_2(-1)=\Pcal'_2(-1)=0$. We put $S_{\tau_i}(x)\coloneqq (1+x)^2R_{\tau_i}(x)$ for $i=1,2,3,4,5$. Note that
\[
\sum_{i=1,2,3,4,5}S_{\tau_i}(x)\equiv 0\bmod{(1+x)^2}
\]
due to Example \ref{exa:the sum of L-functions}. One can see that $H_{\tau_i}(-1)=\prod_\alpha(1+\alpha)^2\eqqcolon A$ for $i=1,2,3,4,5$, and so that
\begin{align*}
\Pcal_2(-1)&=\sum_{i=1,2,3,4,5}S_{\tau_i}(-1)A=0,\\
\Pcal'_2(-1)&=\sum_{i=1,2,3,4,5}(S'_{\tau_i}(-1)A+S_{\tau_i}(-1)H'_{\tau_i}(-1))=\sum_{i=1,2,3,4,5}S_{\tau_i}(-1)H'_{\tau_i}(-1).
\end{align*}
If $J$ contains $-1$, then one sees easily that $H'_{\tau_i}(-1)=0$ for $i=1,2,3,4,5$. Otherwise, one can calculate as
\begin{align*}
&H'_{\tau_1}(-1)=-4AB,\ \ H'_{\tau_1}(-1)=-2AB,\ \ H'_{\tau_3}(-1)=-AB,\\ 
&H'_{\tau_4}(-1)=-AB,\ \ H'_{\tau_5}(-1)=0.
\end{align*}
for $B\coloneqq \sum_\alpha\frac{\alpha}{1+\alpha}$. Therefore,
\begin{align*}
\Pcal'_2(-1)&=H'_{\tau_1}(-1)+H'_{\tau_2}(-1)-4H'_{\tau_3}(-1)-2H'_{\tau_4}(-1)+4H'_{\tau_5}(-1)\\
&=(1\cdot(-4)+1\cdot(-2)-4\cdot(-1)-2\cdot(-1)+4\cdot 0)AB\\
&=0
\end{align*}
and we have checked (b). 

Next, we assume that $q$ is even. We define $R_{\tau'_i}(x)$ as
\begin{align*}
&R_{\tau'_1}(x)=\frac{1+x+x^2}{(1+x)^2(1+x^2)},\ \ R_{\tau'_3}(x)=\frac{x}{(1+x)^2},\\
&R_{\tau'_4}(x)=\frac{x}{(1+x)^2},\ \ R_{\tau'_5}(x)=\frac{x(x-2)}{2(1+x)^2},\ \ R_{\tau'_6}(x)=\frac{x^2}{2(1+x^2)}
\end{align*}
 and $H_{\tau'_i}(x)$ as $H_{\tau_i}(x)$ for $i=1,3,4,5,6$. We put
\[
\Qcal(x)\coloneqq \sum_{i=1,3,4,5,6}R_{\tau'_i}(x)H_{\tau'_i}(x)
\]
so that 
\[
\Qcal(q)=\sum_{[\gamma]}L_S(M_{G_\gamma})=\sum_{i=1,3,4,5,6}N_{\tau'_i}(q)L_S(M_{G_{\tau'_i}}).
\]
Now we are going to check (a)--(c) with each $\Pcal_i(x)$ replaced by $\Qcal_i(x)$. For (a), one has
\begin{align*}
\Qcal_1(x)&\equiv 2(R_{\tau'_5}(x)H_{\tau'_5}(x)+R_{\tau'_6}(x)H_{\tau'_6}(x))\\
&=\frac{x(x-2)}{(1+x)^2}\prod_\alpha(1+\alpha)^2+\frac{x^2}{1+x^2}\prod_\alpha(1+\alpha^2)\\
&\equiv \frac{x^2}{(1+x)^2}\left(\prod_\alpha(1+\alpha)^2+\prod_\alpha(1+\alpha^2)\right)\bmod 2.
\end{align*}
The term inside the parenthesis is zero modulo $2\overline{\Z}$, and thus is zero modulo $2\Z$. 

For (c), one has
\begin{align*}
\Qcal_3(x)&\equiv(1+x^2)(R_{\tau'_1}(x)H_{\tau'_1}(x)+R_{\tau'_6}(x)H_{\tau'_6}(x))\\
&=\frac{1+x+x^2}{(1+x)^2}\prod_\alpha(1-\alpha x)(1-\alpha x^3)+\frac{x^2}{2}\prod_\alpha(1+\alpha^2)\bmod 1+x^2
\end{align*}
and so that 
\begin{align*}
\Qcal_3(\sqrt{-1})=\frac{1}{2}\cdot\prod_\alpha(1+\alpha^2)+\frac{-1}{2}\cdot\prod_\alpha(1+\alpha^2)=0.
\end{align*}

For (b), we put $S_{\tau'_i}(x)\coloneqq (1+x)^2R_{\tau'_i}(x)$ for $i=1,3,4,5$. Then one can calculate as
\begin{align*}
\Qcal_2(-1)&=\sum_{i=1,3,4,5}S_{\tau'_i}(-1)A=0,\\
\Qcal'_2(-1)&=\sum_{i=1,3,4,5}(S'_{\tau'_i}(-1)A+S_{\tau'_i}(-1)H'_{\tau'_i}(-1))=\sum_{i=1,3,4,5}S_{\tau'_i}(-1)H'_{\tau'_i}(-1)\\
&=\left(\frac{1}{2}\cdot(-4)-1\cdot(-1)-1\cdot(-1)+\frac{3}{2}\cdot 0\right)AB=0.
\end{align*}

The above argument shows that $\Pcal(q)$ and $\Qcal(q)$ are polynomials in $\{q\}\cup J$ over $\Z$ and the base change from $\F_q$ to $\F_{q^m}$ replaces $\{q\}\cup J$ by $\{q^m\}\cup J^m$ as in the proof of Theorem \ref{thm:the sum of L-functions for SL_l}. Thus our assertion follows.
\end{prf}

Next we work on the case $G=\mathrm{Sp}_6$. In this case, there are twelve types of semisimple classes if $q$ is odd, and ten types if $q$ is even. (See Table 3, 4 for details.)
\renewcommand{\arraystretch}{1.3}
\begin{table}[ht]
\centering
\caption{types of semisimple classes and related data for $\Sp_6$, $q$ odd}
\resizebox{\textwidth}{!}{
\begin{tabular}{|c|c|c|c|}\hhline{|-|-|-|-|}
    type $\tau$ & $G_\tau$ & $N_\tau(q)$ & $\det(1-t\Fr_q\mid M_{G_\tau})$ \\ \hhline{|=|=|=|=|}
    $\tau_1=(3,0;;)$ & $\Sp_6$ & 2 & $(1-t q)(1-t q^3)(1-t q^5)$\\ \hhline{|-|-|-|-|}
    $\tau_2=(2,1;;)$ & $\Sp_4\times\Sp_2$ & $2$ & $(1-t q)^2(1-t q^3)$\\ \hhline{|-|-|-|-|}
    $\tau_3=(2,0;(1,1);)$ & $\Sp_4\times \U_1$ & $q-1$ & $(1+t)(1-t q)(1-t q^3)$\\ \hhline{|-|-|-|-|}
    $\tau_4=(1,1;(1,1);)$ & $\Sp_2\times\Sp_2\times\U_1$ & $(q-1)/2$ & $(1+t)(1-t q)^2$ \\ \hhline{|-|-|-|-|}
    $\tau_5=(1,0;(1,2);)$ & $\Sp_2\times\U_2$ & $q-1$ & $(1+t)(1-t q)^2$ \\ \hhline{|-|-|-|-|}
    $\tau_6=(1,0;(1,1),(1,1);)$ & $\Sp_2\times\U_1\times\U_1$ & $(q-1)(q-3)/4$ & $(1+t)^2(1-t q)$ \\ \hhline{|-|-|-|-|}
    $\tau_7=(1,0;(2,1);)$ & $\Sp_2\times\Res_{\F_{q^2}/\F_q}\U_1$ & $(q^2-1)/2$ & $(1+t^2)(1-t q)$ \\ \hhline{|-|-|-|-|}
    $\tau_8=(;(1,3);)$ & $\U_3$ & $(q-1)/2$ & $(1+t)(1-t q)(1+t q^2)$ \\ \hhline{|-|-|-|-|}
    $\tau_9=(;(1,2),(1,1);)$ & $\U_2\times\U_1$ & $(q-1)(q-3)/4$ & $(1+t)^2(1-t q)$ \\ \hhline{|-|-|-|-|}
    $\tau_{10}=(;(1,1),(1,1),(1,1);)$ & $\U_1\times\U_1\times\U_1$ & $(q-1)(q-3)(q-5)/48$ & $(1+t)^3$ \\ \hhline{|-|-|-|-|}
    $\tau_{11}=(;(2,1),(1,1);)$ & $\Res_{\F_{q^2}/\F_q}\U_1\times\U_1$ & $(q-1)(q^2-1)/8$ & $(1+t)(1+t^2)$ \\ \hhline{|-|-|-|-|}
    $\tau_{12}=(;(3,1);)$ & $\Res_{\F_{q^3}/\F_q}\U_1$ & $(q^3-q)/6$ & $1+t^3$ \\ \hhline{|-|-|-|-|}
\end{tabular}}

\caption{types of semisimple classes and related data for $\Sp_6$, $q$ even}
\resizebox{\textwidth}{!}{
\begin{tabular}{|c|c|c|c|}\hhline{|-|-|-|-|}
    type $\tau'$ & $G_{\tau'}$ & $N_{\tau'}(q)$ & $\det(1-t\Fr_q\mid M_{G_{\tau'}})$ \\ \hhline{|=|=|=|=|}
    $\tau'_1=(3;;)$ & $\Sp_6$ & 1 & $(1-t q)(1-t q^3)(1-t q^5)$\\ \hhline{|-|-|-|-|}
    $\tau'_3=(2;(1,1);)$ & $\Sp_4\times \U_1$ & $q/2$ & $(1+t)(1-t q)(1-t q^3)$\\ \hhline{|-|-|-|-|}
    $\tau'_5=(1;(1,2);)$ & $\Sp_2\times\U_2$ & $q/2$ & $(1+t)(1-t q)^2$ \\ \hhline{|-|-|-|-|}
    $\tau'_6=(1;(1,1),(1,1);)$ & $\Sp_2\times\U_1\times\U_1$ & $q(q-2)/8$ & $(1+t)^2(1-t q)$ \\ \hhline{|-|-|-|-|}
    $\tau'_7=(1;(2,1);)$ & $\Sp_2\times\Res_{\F_{q^2}/\F_q}\U_1$ & $q^2/4$ & $(1+t^2)(1-t q)$ \\ \hhline{|-|-|-|-|}
    $\tau'_8=(;(1,3);)$ & $\U_3$ & $q/2$ & $(1+t)(1-t q)(1+t q^2)$ \\ \hhline{|-|-|-|-|}
    $\tau'_9=(;(1,2),(1,1);)$ & $\U_2\times\U_1$ & $q(q-2)/4$ & $(1+t)^2(1-t q)$ \\ \hhline{|-|-|-|-|}
    $\tau'_{10}=(;(1,1),(1,1),(1,1);)$ & $\U_1\times\U_1\times\U_1$ & $q(q-2)(q-4)/48$ & $(1+t)^3$ \\ \hhline{|-|-|-|-|}
    $\tau'_{11}=(;(2,1),(1,1);)$ & $\Res_{\F_{q^2}/\F_q}\U_1\times\U_1$ & $q^3/8$ & $(1+t)(1+t^2)$ \\ \hhline{|-|-|-|-|}
    $\tau'_{12}=(;(3,1);)$ & $\Res_{\F_{q^3}/\F_q}\U_1$ & $(q^3-q)/6$ & $1+t^3$ \\ \hhline{|-|-|-|-|}
\end{tabular}}
\end{table}
\renewcommand{\arraystretch}{1}

\begin{thm}\label{thm:the sum of L-functions for Sp_6}
The function $m\mapsto\mathcal{L}(\mathrm{Sp}_6,m)$ is of Lefschetz type. 
\end{thm}

\begin{prf}
First, we assume that $q$ is odd. We define the rational functions $R_{\tau_i}(x)\in\Z(x)$ as
\begin{align*}
&R_{\tau_1}(x)=\frac{2(1+x+x^2+x^3+x^4)}{(1+x)^3(1+x^2)(1-x+x^2)},\ \ R_{\tau_2}(x)=\frac{2(1+x+x^2)}{(1+x)^3(1+x^2)},\\
&R_{\tau_3}(x)=\frac{2(x-1)(1+x+x^2)}{(1+x)^3(1+x^2)},\ \ R_{\tau_4}(x)=\frac{x-1}{(1+x)^3},\ \ R_{\tau_5}(x)=\frac{2(x-1)}{(1+x)^3},\\
&R_{\tau_6}(x)=\frac{(x-1)(x-3)}{(1+x)^3},\ \ R_{\tau_7}(x)=\frac{x-1}{1+x^2},\ \ R_{\tau_8}(x)=\frac{(x-1)(1+x^2)}{(1+x)^3(1-x+x^2)}\\
&R_{\tau_9}(x)=\frac{(x-1)(x-3)}{(1+x)^3},\ \ R_{\tau_{10}}(x)=\frac{(x-1)(x-3)(x-5)}{6(1+x)^3},\\
&R_{\tau_{11}}(x)=\frac{(x-1)^2}{2(1+x^2)},\ \ R_{\tau_{12}}(x)=\frac{x(x-1)}{3(1-x+x^2)}
\end{align*}
and the polynomials $H_{\tau_i}(x)\in\Z[x]$ as
\begin{align*}
&H_{\tau_1}(x)=\prod_{\alpha\in J}(1-\alpha x)(1-\alpha x^3)(1-\alpha x^5),\ \ H_{\tau_2}(x)=\prod_{\alpha\in J}(1-\alpha x)^2(1-\alpha x^3),\\
&H_{\tau_3}(x)=\prod_{\alpha\in J}(1+\alpha)(1-\alpha x)(1-\alpha x^3),\ \ H_{\tau_4}(x)=\prod_{\alpha\in J}(1+\alpha)(1-\alpha x)^2,\\
&H_{\tau_5}(x)=\prod_{\alpha\in J}(1+\alpha)(1-\alpha x)^2,\ \ H_{\tau_6}(x)=\prod_{\alpha\in J}(1+\alpha)^2(1-\alpha x),\\
&H_{\tau_7}(x)=\prod_{\alpha\in J}(1+\alpha^2)(1-\alpha x),\ \ H_{\tau_8}(x)=\prod_{\alpha\in J}(1+\alpha)(1-\alpha x)(1+\alpha x^2),\\
&H_{\tau_9}(x)=\prod_{\alpha\in J}(1+\alpha)^2(1-\alpha x),\ \ H_{\tau_{10}}(x)=\prod_{\alpha\in J}(1+\alpha)^3,\\
&H_{\tau_{11}}(x)=\prod_{\alpha\in J}(1+\alpha)(1+\alpha^2),\ \ H_{\tau_{12}}(x)=\prod_{\alpha\in J}(1+\alpha^3).
\end{align*}
 Then one can see that 
\begin{align*}
&R_{\tau_i}(q)= N_{\tau_i}(q)\frac{\det(1-\Fr_q\mid M_{G_{\tau_i}})}{\det(1-q\Fr_q\mid M_{G_{\tau_i}})},\\
&H_{\tau_i}(q)= \prod_{\alpha\in J}\det(1-\alpha \Fr_q\mid M_{G_{\tau_i}})
\end{align*}
for $i=1,\ldots, 12$ using Table 3. We put
\[
\Pcal(x)\coloneqq \sum_{i=1}^{12} R_{\tau_i}(x)H_{\tau_i}(x)\in\Z(x)
\]
so that
\[
\Pcal(q)=\sum_{[\gamma]}L_S(M_{G_\gamma})=\sum_{i=1}^{12}N_{\tau_i}(q)L_S(M_{G_{\tau_i}}).
\]
Now we are going to prove that $\Pcal(x)\in\Z[x]$. One can see that $\Pcal(x)$ lies in $(6(1+x)^3(1+x^2)(1-x+x^2))^{-1}\Z[x]$. we put $\Pcal_1(x)\coloneqq 2\Pcal(x)$, $\Pcal_2(x)\coloneqq 3\Pcal(x)$, $\Pcal_3(x)\coloneqq (1+x)^3\Pcal(x)$, $\Pcal_4(x)\coloneqq (1+x^2)\Pcal(x)$, $\Pcal_5(x)\coloneqq (1-x+x^2)\Pcal(x)$ and show that
\begin{enumerate}
\item[(a)] $\Pcal_1(x)\equiv 0\bmod 2$ in $(3(1+x)^3(1+x^2)(1-x+x^2))^{-1}\Z[x]$,
\item[(b)] $\Pcal_2(x)\equiv 0\bmod 3$ in $(2(1+x)^3(1+x^2)(1-x+x^2))^{-1}\Z[x]$,
\item[(c)] $\Pcal_3(x)\equiv 0\bmod (1+x)^3$ in $(6(1+x^2)(1-x+x^2))^{-1}\Z[x]$,
\item[(c)] $\Pcal_4(x)\equiv 0\bmod 1+x^2$ in $(6(1+x)^3(1-x+x^2))^{-1}\Z[x]$,
\item[(e)] $\Pcal_4(x)\equiv 0\bmod 1-x+x^2$ in $(6(1+x)^3(1+x^2))^{-1}\Z[x]$.
\end{enumerate}

For (a), one has
\begin{align*}
\Pcal_1(x)&\equiv 2(R_{\tau_{10}}(x)H_{\tau_{10}}(x)+R_{\tau_{11}}(x)H_{\tau_{11}}(x))\\
&=\frac{(x-1)(x-3)(x-5)}{3(1+x)^3}\prod_{\alpha}(1+\alpha)^3+\frac{(x-1)^2}{1+x^2}\prod_\alpha(1+\alpha)(1+\alpha^2)\\
&\equiv \prod_{\alpha}(1+\alpha)^3+\prod_\alpha(1+\alpha)(1+\alpha^2)\bmod 2.
\end{align*}
It is zero modulo $2\overline{\Z}$, and thus is zero modulo $2\Z$. 

For (b), one has
\begin{align*}
\Pcal_2(x)&\equiv 3(R_{\tau_{10}}(x)H_{\tau_{10}}(x)+R_{\tau_{12}}(x)H_{\tau_{12}}(x))\\
&=\frac{(x-1)(x-3)(x-5)}{2(1+x)^3}\prod_{\alpha}(1+\alpha)^3+\frac{x(x-1)}{1-x+x^2}\prod_\alpha(1+\alpha^3)\\
&\equiv \frac{x(x-1)}{(1+x)^2}\left(-\prod_{\alpha}(1+\alpha)^3+\prod_\alpha(1+\alpha^3)\right)\bmod 3.
\end{align*}
The term inside the parenthesis is zero modulo $3\overline{\Z}$, and thus is zero modulo $3\Z$. 

For (d), it suffices to show that $\Pcal_4(\sqrt{-1})=0$. One has
\begin{align*}
\Pcal_4(x)&\equiv(1+x^2)\sum_{i=1,2,3,7,11}R_{\tau_i}(x)H_{\tau_i}(x)\bmod 1+x^2
\end{align*}
and so that
\begin{align*}
\Pcal_4(\sqrt{-1})&=\left(\frac{1-\sqrt{-1}}{2}+\frac{1-\sqrt{-1}}{2}+(\sqrt{-1}-1)\right)\cdot\prod_{\alpha}(1-\sqrt{-1}\alpha)(1+\alpha^2)\\
&\quad\quad+\left(\sqrt{-1}-\sqrt{-1}\right)\cdot\prod_\alpha(1+\alpha)(1+\alpha^2)\\
&=0.
\end{align*}

For (e), it suffices to show that $\Pcal_5(\zeta_6)=0$ for $\zeta_6=e^{2\pi\sqrt{-1}/6}$. One has
\begin{align*}
\Pcal_5(x)&\equiv(1-x+x^2)\sum_{i=1,8,12}R_{\tau_i}(x)H_{\tau_i}(x)\bmod 1-x+x^2
\end{align*}
and so that
\begin{align*}
\Pcal_5(\zeta_6)=\left(\left(\frac{1}{3}-\frac{\sqrt{-3}}{9}\right)+\frac{\sqrt{-3}}{9}-\frac{1}{3}\right)\prod_\alpha(1+\alpha^3)=0.
\end{align*}

For (c), it suffices to show that $\Pcal_3(-1)=\Pcal'_3(-1)=\Pcal''_3(-1)=0$. Put $S_{\tau_i}(x)\coloneqq (1+x)^3R_{\tau_i}(x)$ for $i=1,2,3,4,5,6,8,9,10$. Note that
\[
\sum_{i=1,2,3,4,5,6,8,9,10}S_{\tau_i}(x)\equiv 0\bmod (1+x)^3
\]
due to Example \ref{exa:the sum of L-functions}. One can see that $H_{\tau_i}(-1)=\prod_\alpha(1+\alpha)^3\eqqcolon A$ for $i=1,2,3,4,5,6,8,9,10$, and so that
\begin{align*}
\Pcal_3(-1)&=\sum_{i=1,2,3,4,5,6,8,9,10}S_{\tau_i}(-1)A=0,\\
\Pcal'_3(-1)&=\sum_{i=1,2,3,4,5,6,8,9,10}(S'_{\tau_i}(-1)A+S_{\tau_i}(-1)H'_{\tau_i}(-1))\\
&=\sum_{i=1,2,3,4,5,6,8,9,10}S_{\tau_i}(-1)H'_{\tau_i}(-1),\\
\Pcal''_3(-1)&=\sum_{i=1,2,3,4,5,6,8,9,10}(S''_{\tau_i}(-1)A+2 S'_{\tau_i}(-1)H'_{\tau_i}(-1)+S_{\tau_i}(-1)H''_{\tau_i}(-1))\\
&=2\sum_{i=1,2,3,4,5,6,8,9,10}S'_{\tau_i}(-1)H'_{\tau_i}(-1)+\sum_{i=1,2,3,4,5,6,8,9,10}S_{\tau_i}(-1)H''_{\tau_i}(-1).
\end{align*}
If $J$ contains $-1$, then one sees easily that $H'_{\tau_i}(-1)=H''_{\tau_i}(-1)=0$ for $i=1,2,3,4,5,6,8,9,10$. Otherwise, one can calculate as
\renewcommand{\arraystretch}{1.3}
\begin{center}
\begin{tabular}{|c|c|c|c|c|}\hhline{|-|-|-|-|-|}
	$\tau$ & $S_{\tau}(-1)$ & $S'_{\tau}(-1)$ & $H'_{\tau}(-1)$ & $H''_{\tau}(-1)$\\ \hhline{|=|=|=|=|=|}
	$\tau_1$ & $1/3$ & $0$ & $-9AB$ & $26AB+46AC+81AD$ \\ \hhline{|-|-|-|-|-|}
	$\tau_2$ & $1$ & $0$ & $-5AB$ & $6AB+14AC+25AD$ \\ \hhline{|-|-|-|-|-|}
	$\tau_3$ & $-2$ & $1$ & $-4AB$ & $6AB+6AC+16AD$ \\ \hhline{|-|-|-|-|-|}
	$\tau_4$ & $-2$ & $1$ & $-2AB$ & $2AC+4AD$ \\ \hhline{|-|-|-|-|-|}
	$\tau_5$ & $-4$ & $2$ & $-2AB$ & $2AC+4AD$ \\ \hhline{|-|-|-|-|-|}
	$\tau_6$ & $8$ & $-6$ & $-AB$ & $AD$ \\ \hhline{|-|-|-|-|-|}
	$\tau_8$ & $-4/3$ & $2/3$ & $-3AB$ & $2AB+4AC+9AD$ \\ \hhline{|-|-|-|-|-|}
	$\tau_9$ & $8$ & $-6$ & $-AB$ & $AD$ \\ \hhline{|-|-|-|-|-|}
	$\tau_{10}$ & $-8$ & $22/3$ & $0$ & $0$ \\ \hhline{|-|-|-|-|-|}
\end{tabular}
\end{center}
\renewcommand{\arraystretch}{1}
where 
\[
B\coloneqq\sum_\alpha\frac{\alpha}{1+\alpha},\ \ C\coloneqq\sum_\alpha\left(\frac{\alpha}{1+\alpha}\right)^2,\ \ D\coloneqq\sum_{\alpha,\beta,\alpha\neq\beta}\frac{\alpha}{1+\alpha}\frac{\beta}{1+\beta}.
\]
Note that $\alpha\neq\beta$ means that $\alpha$ and $\beta$ are different elements of the multiset $J$, and they might be the same element of $\overline{\Z}$. Now we can calculate as
\begin{align*}
\Pcal'_3(-1)&=\left(\frac{1}{3}\cdot(-9)+1\cdot(-5)+(-2)\cdot(-4)+(-2)\cdot(-2)+(-4)\cdot(-2)\right.\\
&\quad\quad\left.{}+8\cdot(-1)+\left(-\frac{4}{3}\right)\cdot(-3)+8\cdot(-1)+(-8)\cdot0\right)AB\\
&=0,\\
\Pcal''_3(-1)&=2\bigg(0\cdot(-9)+0\cdot(-5)+1\cdot(-4)+1\cdot(-2)+2\cdot(-2)\\
&\quad\quad\left.{}+(-6)\cdot(-1)+\frac{2}{3}\cdot(-3)+(-6)\cdot(-1)+\frac{22}{3}\cdot0\right)AB\\
&\quad+\left(\frac{1}{3}\cdot26+1\cdot6+(-2)\cdot6+(-2)\cdot0+(-4)\cdot0\right.\\
&\quad\quad\left.{}+8\cdot0+\left(-\frac{4}{3}\right)\cdot2+8\cdot0+(-8)\cdot0\right)AB\\
&\quad+\left(\frac{1}{3}\cdot46+1\cdot14+(-2)\cdot6+(-2)\cdot2+(-4)\cdot2\right.\\
&\quad\quad\left.{}+8\cdot0+\left(-\frac{4}{3}\right)\cdot4+8\cdot0+(-8)\cdot0\right)AC\\
&\quad+\left(\frac{1}{3}\cdot81+1\cdot25+(-2)\cdot16+(-2)\cdot4+(-4)\cdot4\right.\\
&\quad\quad\left.{}+8\cdot1+\left(-\frac{4}{3}\right)\cdot9+8\cdot1+(-8)\cdot0\right)AD\\
&=0
\end{align*}
and we have checked (c). 

Next, we assume that $q$ is even. We define $R_{\tau'_i}(x)$ as
\begin{align*}
&R_{\tau'_1}(x)=\frac{1+x+x^2+x^3+x^4}{(1+x)^3(1+x^2)(1-x+x^2)},\ \ R_{\tau'_3}(x)=\frac{x(1+x+x^2)}{(1+x)^3(1+x^2)},\\
&R_{\tau'_5}(x)=\frac{x}{(1+x)^3},\ \ R_{\tau'_6}(x)=\frac{x(x-2)}{2(1+x)^3},\ \ R_{\tau'_7}(x)=\frac{x^2}{2(1+x)(1+x^2)},\\
&R_{\tau'_8}(x)=\frac{x(1+x^2)}{(1+x)^3(1-x+x^2)},\ \ R_{\tau'_9}(x)=\frac{x(x-2)}{(1+x)^3},\\
&R_{\tau'_{10}}(x)=\frac{x(x-2)(x-4)}{6(1+x)^3},\ \ R_{\tau'_{11}}(x)=\frac{x^3}{2(1+x)(1+x^2)},\\
&R_{\tau'_{12}}(x)=\frac{x(x-1)}{3(1-x+x^2)}
\end{align*}
 and $H_{\tau'_i}(x)$ as $H_{\tau_i}(x)$ for $i=1,3,4,5,6,7,8,9,10,11,12$. We put 
\[
\Qcal(x)\coloneqq \sum_{i=1,3,5,6,7,8,9,10,11,12}R_{\tau'_i}(x)H_{\tau'_i}(x)
\]
so that
\[
\Qcal(q)=\sum_{[\gamma]}L_S(M_{G_\gamma})=\sum_{i=1,3,5,6,7,8,9,10,11,12}N_{\tau_i}(q)L_S(M_{G_{\tau_i}}).
\]
Now we are going to check (a)--(e) with each $\Pcal_i(x)$ replaced by $\Qcal_i(x)$. For (a), one has
\begin{align*}
\Qcal_1(x)&\equiv 2\sum_{i=6,7,10,11}R_{\tau'_i}(x)H_{\tau'_i}(x)\\
&=\left(\frac{x(x-2)}{(1+x)^3}+\frac{x^2}{(1+x)(1+x^2)}\right)\prod_\alpha(1+\alpha)^2(1-\alpha x)\\
&\quad +\frac{x(x-2)(x-4)}{3(1+x)^3}\prod_{\alpha}(1+\alpha)^3+\frac{x^3}{(1+x)(1+x^2)}\prod_\alpha(1+\alpha)(1+\alpha^2)\\
&\equiv \frac{x^3}{(1+x)^3}\left(\prod_{\alpha}(1+\alpha)^3+\prod_\alpha(1+\alpha)(1+\alpha^2)\right)\bmod 2.
\end{align*}
The term inside the parenthesis is zero modulo $2\overline{\Z}$, and thus is zero modulo $2\Z$. 

For (b), one has
\begin{align*}
\Qcal_2(x)&\equiv 3(R_{\tau'_{10}}(x)H_{\tau'_{10}}(x)+R_{\tau'_{12}}(x)H_{\tau'_{12}}(x))\\
&=\frac{x(x-2)(x-4)}{2(1+x)^3}\prod_{\alpha}(1+\alpha)^3+\frac{x(x-1)}{1-x+x^2}\prod_\alpha(1+\alpha^3)\\
&\equiv \frac{x(x-1)}{(1+x)^2}\left(-\prod_{\alpha}(1+\alpha)^3+\prod_\alpha(1+\alpha^3)\right)\bmod 3.
\end{align*}
The term inside the parenthesis is zero modulo $3\overline{\Z}$, and thus is zero modulo $3\Z$. 

For (d), one has
\begin{align*}
\Qcal_4(x)&\equiv(1+x^2)\sum_{i=1,3,7,11}R_{\tau_i}(x)H_{\tau_i}(x)\bmod 1+x^2
\end{align*}
and so that
\begin{align*}
\Qcal_4(\sqrt{-1})&=\left(\frac{1-\sqrt{-1}}{4}+\frac{-1+\sqrt{-1}}{4}\right)\cdot\prod_\alpha(1-\sqrt{-1\alpha})(1+\alpha^2)\\
&\quad\quad+\left(\frac{1+\sqrt{-1}}{4}+\frac{-1-\sqrt{-1}}{4}\right)\cdot\prod_\alpha(1+\alpha)(1+\alpha^2)\\
&=0.
\end{align*}

For (e), one has
\begin{align*}
\Qcal_5(x)&\equiv(1-x+x^2)\sum_{i=1,8,12}R_{\tau_i}(x)H_{\tau_i}(x)\bmod 1-x+x^2
\end{align*}
and so that
\begin{align*}
\Qcal_5(\zeta_6)=\left(\left(\frac{1}{6}-\frac{\sqrt{-3}}{18}\right)+\left(\frac{1}{6}+\frac{\sqrt{-3}}{18}\right)-\frac{1}{3}\right)\prod_\alpha(1+\alpha^3)=0.
\end{align*}

For (c), we put $S_{\tau'_i}(x)\coloneqq (1+x)^3R_{\tau'_i}(x)$ for $i=1,3,5,6,7,8,9,10,11$ and $A'\coloneqq\prod_\alpha(1+\alpha)(1+\alpha^2)$. Note that $H_{\tau'_i}(-1)=A'$ and $S_{\tau'_i}(-1)=S'_{\tau'_i}(-1)=0$ for $i=7,11$. Then one can calculate as
\begin{align*}
\Qcal_3(-1)&=\sum_{i=1,3,5,6,8,9,10}S_{\tau'_i}(-1)A=0,\\
\Qcal'_3(-1)&=\sum_{i=1,3,5,6,8,9,10}(S'_{\tau'_i}(-1)A+S_{\tau'_i}(-1)H'_{\tau'_i}(-1))\\
&=\sum_{i=1,3,5,6,8,9,10}S_{\tau'_i}(-1)H'_{\tau'_i}(-1),\\
\Qcal''_3(-1)&=\sum_{i=1,3,5,6,8,9,10}(S''_{\tau'_i}(-1)A+2S'_{\tau'_i}(-1)H'_{\tau'_i}(-1)+S_{\tau'_i}(-1)H''_{\tau'_i}(-1))\\
&\quad\quad+\sum_{i=7,11}S''_{\tau'_i}(-1)A'\\
&=2\sum_{i=1,3,5,6,8,9,10}S'_{\tau'_i}(-1)H'_{\tau'_i}(-1)+\sum_{i=1,3,5,6,8,9,10}S_{\tau'_i}(-1)H''_{\tau'_i}(-1)\\
&\quad\quad+\sum_{i=7,11}S''_{\tau'_i}(-1)(A'-A).
\end{align*}
If $J$ contains $-1$, then one sees easily that $H'_{\tau_i}(-1)=H''_{\tau_i}(-1)=0$ for $i=1,3,5,6,8,9,10$. Otherwise, one can calculate as
\renewcommand{\arraystretch}{1.3}
\begin{center}
\begin{tabular}{|c|c|c|c|c|}\hhline{|-|-|-|-|-|}
	$\tau'$ & $S_{\tau'}(-1)$ & $S'_{\tau'}(-1)$ & $H'_{\tau'}(-1)$ & $H''_{\tau'}(-1)$\\ \hhline{|=|=|=|=|=|}
	$\tau'_1$ & $1/6$ & $0$ & $-9AB$ & $26AB+46AC+81AD$ \\ \hhline{|-|-|-|-|-|}
	$\tau'_3$ & $-1/2$ & $1/2$ & $-4AB$ & $6AB+6AC+16AD$ \\ \hhline{|-|-|-|-|-|}
	$\tau'_5$ & $-1$ & $1$ & $-2AB$ & $2AC+4AD$ \\ \hhline{|-|-|-|-|-|}
	$\tau'_6$ & $3/2$ & $-2$ & $-AB$ & $AD$ \\ \hhline{|-|-|-|-|-|}
	$\tau'_8$ & $-2/3$ & $2/3$ & $-3AB$ & $2AB+4AC+9AD$ \\ \hhline{|-|-|-|-|-|}
	$\tau'_9$ & $3$ & $-4$ & $-AB$ & $AD$ \\ \hhline{|-|-|-|-|-|}
	$\tau'_{10}$ & $-5/2$ & $23/6$ & $0$ & $0$ \\ \hhline{|-|-|-|-|-|}
\end{tabular}
\end{center}
\renewcommand{\arraystretch}{1}
and $S''_{\tau'_7}(-1)=1/2$, $S''_{\tau'_{11}}(-1)=-1/2$. It follows that
\begin{align*}
\Qcal'_3(-1)&=\left(\frac{1}{6}\cdot(-9)+\left(-\frac{1}{2}\right)\cdot(-4)+(-1)\cdot(-2)+\frac{3}{2}\cdot(-1)\right.\\
&\quad\quad\left.+\left(-\frac{2}{3}\right)\cdot(-3)+3\cdot(-1)+\left(-\frac{5}{2}\right)\cdot0\right)AB\\
&=0,\\
\Qcal''_3(-1)&=2\left(0\cdot(-9)+\frac{1}{2}\cdot(-4)+1\cdot(-2)+(-2)\cdot(-1)\right.\\
&\quad\quad\left.{}+\frac{2}{3}\cdot(-3)+(-4)\cdot(-1)+\frac{23}{6}\cdot0\right)AB\\
&\quad+\left(\frac{1}{6}\cdot26+\left(-\frac{1}{2}\right)\cdot6+(-1)\cdot0+\frac{3}{2}\cdot0\right.\\
&\quad\quad\left.+\left(-\frac{2}{3}\right)\cdot2+3\cdot0+\left(-\frac{5}{2}\right)\cdot0\right)AB\\
&\quad+\left(\frac{1}{6}\cdot46+\left(-\frac{1}{2}\right)\cdot6+(-1)\cdot2+\frac{3}{2}\cdot0\right.\\
&\quad\quad\left.+\left(-\frac{2}{3}\right)\cdot4+3\cdot0+\left(-\frac{5}{2}\right)\cdot0\right)AC\\
&\quad+\left(\frac{1}{6}\cdot81+\left(-\frac{1}{2}\right)\cdot16+(-1)\cdot4+\frac{3}{2}\cdot1\right.\\
&\quad\quad\left.+\left(-\frac{2}{3}\right)\cdot9+3\cdot1+\left(-\frac{5}{2}\right)\cdot0\right)AD\\
&=0
\end{align*}

The above argument shows that $\Pcal(q)$ and $\Qcal(q)$ are polynomials in $\{q\}\cup J$ over $\Z$ and the base change from $\F_q$ to $\F_{q^m}$ replaces $\{q\}\cup J$ by $\{q^m\}\cup J^m$ . Thus our assertion follows.
\end{prf}

\begin{rem}
We verified that Conjecture \ref{conj:the sum of L-functions} is true for $G=\Sp_4,\Sp_6$. It follows form Theorem \ref{thm:the sum of multiplicities} that the functions
\[
m\mapsto \sum_{\pi\in\mathrm{Irr}_0^{S_m}(\Sp_{4,k_m})}m(\pi),\ \ \ m\mapsto \sum_{\pi\in\mathrm{Irr}_0^{S_m}(\Sp_{6,k_m})}m(\pi)
\]
are of Lefschetz type if we assume Conjectures \ref{conj:EP function} and \ref{conj:trace formula}.

It is expected that Conjecture \ref{conj:the sum of L-functions} is true for $G=\Sp_{2n}$ in general, and that there is a comprehensive proof for it. We have that there exists a rational function $\Pcal(x,J)\in\Z(x,J)$ such that $\mathcal{L}(\Sp_{2n},m)=\Pcal(q^m,J^m)$ as in the proofs of Theorems \ref{thm:the sum of L-functions for Sp_4} and \ref{thm:the sum of L-functions for Sp_6}. Therefore, it is enough to show that $\Pcal(x,J)$ lies in $\Z[x,J]$.
\end{rem}

\bibliographystyle{my_amsalpha}
\bibliography{master}

\end{document}